\documentclass{amsart}
\usepackage{amssymb,amsxtra,latexsym,graphics}

\newcommand{\F}{\mathcal{F}}
\newcommand{\G}{\mathcal{G}}
\newcommand{\A}{\mathcal{A}}
\newcommand{\B}{\mathcal{B}}
\renewcommand{\L}{\mathcal{L}}
\newcommand{\M}{\mathcal{M}}
\newcommand{\C}{\mathcal{C}}
\newcommand{\D}{\mathcal{D}}
\newcommand{\R}{\mathbb{R}}
\newcommand{\CC}{\mathbb{C}}

\newcommand{\fq}{\mathbb{F}_q}
\newcommand{\abs}[1]{\lvert#1\rvert}
\newcommand{\norm}[1]{\lVert#1\rVert}
\newcommand{\zerohat}{\hat{0}}
\newcommand{\onehat}{\hat{1}}
\newcommand{\supp}{\operatorname{supp}}
\newcommand{\rank}{\operatorname{rank}}

\newcommand{\gl}{\operatorname{GL}}
\newcommand{\des}{\operatorname{des}}
\newcommand{\type}{\operatorname{type}}
\newcommand{\End}{\operatorname{End}}
\newcommand{\inv}{\operatorname{inv}}
\newcommand{\onto}{\twoheadrightarrow}

\newcommand{\iso}{\cong}
\newcommand{\isoto}{\stackrel{\cong}{\longrightarrow}}
\newcommand{\J}{\mathcal{J}}
\newcommand{\qbinom}[2]{\genfrac{[}{]}{0pt}{}{#1}{#2}}

\newcommand{\lhor}{\rule[.7ex]{.2in}{.4pt}}

\theoremstyle{plain}
\newtheorem*{corollary}{Corollary}

\newtheorem{theorem}{Theorem}

\newtheorem*{main}{Main theorem}
\newtheorem{proposition}{Proposition}

\theoremstyle{definition}

\theoremstyle{remark}
\newtheorem*{remark}{Remark}
\newtheorem{example}{Example}

\newcommand{\abels}{abels91}
\newcommand{\riffle}{bayerdiaconis92:_trail}
\newcommand{\bid}{bidigare97:_hyper}
\newcommand{\bergeron}{bergeron92:_descent}
\newcommand{\bhr}{bidigarehanlon98}
\newcommand{\bbd}{billera99:_random}
\newcommand{\bl}{billera98:_noncom}
\newcommand{\red}{bjoerner93:_om}
\newcommand{\brown}{brown89:_build}
\newcommand{\bd}{browndiaconis:_random}

\newcommand{\desar}{desarmenien83:_une}
\newcommand{\desarw}{desarmenien93:_descen}
\newcommand{\drep}{diaconis88:_group}
\newcommand{\dicm}{diaconis98:icm}

\newcommand{\dfp}{diaconisfill92:_top_to_random}
\newcommand{\diacfreed}{diaconis:_iterat}

\newcommand{\Fill}{fill96}
\newcommand{\greene}{greene73:_moebius}
\newcommand{\grillet}{grillet95:_semig}
\newcommand{\gb}{grove85:_finit}
\newcommand{\mukherjea}{hoegnaes95:_probab}
\newcommand{\humphreys}{humphreys90:_reflec_coxet}

\newcommand{\klein}{klein-barmen40:_uber_veral_verban}
\newcommand{\ot}{orlikterao92:_book}
\newcommand{\petrich}{petrich71}
\newcommand{\petrichbook}{petrich77:_lectur}
\newcommand{\phatarfod}{phatarfod91:_matrix}

\newcommand{\schutz}{schuetzenberger47:_sur}
\newcommand{\serre}{serre77:_linear}
\newcommand{\solomonb}{solomon67:_burns}
\newcommand{\solomon}{solomon76:_mackey}
\newcommand{\stanleyflag}{stanley96:_combin}
\newcommand{\stanley}{stanley97:_enumer1}
\newcommand{\tits}{tits74:_build_bn}
\newcommand{\wachs}{wachs89}
\newcommand{\welsh}{welsh76:_matroid}
\newcommand{\whitney}{whitney35}
\newcommand{\zie}{ziegler95:_lectures}
\newcommand{\zas}{zaslavsky75:_facing}
\newcommand{\zasgen}{zaslavsky77}

\begin{document}
\title{Semigroups, rings, and Markov chains}
\author{Kenneth S. Brown}
\address{Department of Mathematics\\
Cornell University\\
Ithaca, NY 14853}
\email{kbrown@math.cornell.edu}

\date{September 20, 1999}

\begin{abstract}
We analyze random walks on a class of semigroups called ``left-regular
bands''.  These walks include the hyperplane chamber walks of
Bidigare, Hanlon, and Rockmore.  Using methods of ring theory, we show
that the transition matrices are diagonalizable and we calculate the
eigenvalues and multiplicities.  The methods lead to explicit formulas
for the projections onto the eigenspaces.  As examples of these
semigroup walks, we construct a random walk on the maximal chains of
any distributive lattice, as well as two random walks associated with
any matroid.  The examples include a $q$-analogue of the Tsetlin
library.  The multiplicities of the eigenvalues in the matroid walks
are ``generalized derangement numbers'', which may be of independent
interest.
\end{abstract}

\maketitle

\section{Introduction} \label{s:intro}

There are many tools available for the study of random walks on finite
groups, an important one being representation theory~\cite{\drep}.
For finite semigroups, on the other hand, there is no representation
theory comparable to that for groups.  And, although there is some
general theory of random walks~\cite{\mukherjea,\diacfreed}, much less
is known for semigroups than for groups.  We consider here a special
class of finite semigroups whose irreducible representations can be
worked out explcitly (they are all 1-dimensional), and we use this
information to analyze the random walks.  In particular, we calculate
the eigenvalues, which turn out to be real.

The semigroups we treat are called ``left-regular bands'' in the
semigroup literature.  There are many interesting examples of them,
including the hyperplane chamber walks introduced by Bidigare, Hanlon,
and Rockmore~\cite{\bhr}, as well as several new examples.  Our approach
via representation theory provides a clear conceptual explanation for
some of the remarkable features of the hyperplane chamber walks proved
in \cite{\bhr,\bd}.

\subsection{Random walks on left-regular bands} \label{sub:semi}

A \emph{left-regular band}, or LRB, is a semigroup $S$ that satisfies
the identities
\begin{equation} \tag{D}
x^2=x \quad \text{and} \quad xyx=xy
\end{equation}
for all $x,y\in S$.  We call (D) the ``deletion property'', because it
admits the following restatement:  Whenever we have a product $x_1
x_2\cdots x_n$ in~$S$, we can delete any factor that has occurred
earlier without changing the value of the product.  Information about
LRBs can be found in \cite{\grillet,\petrich,\petrichbook}.  Early
references to the identity $xyx=xy$ are \cite{\klein,\schutz}.

Our LRBs will always be finite and, for simplicity, will usually have
an identity.  The second assumption involves no loss of generality,
since we can always adjoin an identity to $S$ and property~(D) still
holds.  And even the first assumption involves very little loss of
generality, since (D) implies that $S$ is finite if it is finitely
generated.

To run a random walk on $S$, start with a probability distribution
$\{w_x\}_{x\in S}$ on~$S$.  A step in the walk then goes from $s$
to~$xs$, where $x\in S$ is chosen with probability~$w_x$.  Thus there
is a transition from $s$ to~$t$ with probability
\begin{equation} \label{e:transition1}
P(s,t)= \sum_{xs=t} w_x.
\end{equation}
As we will see in the examples below, it is natural to consider a
slight variant of this walk, in which we confine ourselves to elements
of a left ideal $I\subseteq S$, i.e., a nonempty subset that is
closed under left-multiplication by arbitrary elements of~$S$.  If the
walk starts in~$I$ then it stays there, so we have a Markov chain
on~$I$ with transition matrix given by~\eqref{e:transition1} for
$s,t\in I$.

The next three subsections give examples of LRBs and the associated
random walks.

\subsection{Example: Hyperplane face semigroups} \label{sub:bhr}

These are the motivating examples that led to the present paper.
Briefly, a finite set of affine hyperplanes in a real vector space $V$
divides $V$ into regions called \emph{chambers}.  These are polyhedral
sets, which have faces.  The totality $\F$ of all the faces is a poset
under the face relation.  Less obviously, $\F$ admits a product,
making it a LRB.  See Appendix~\ref{app:hyperplane} for details.
Assume for simplicity that the arrangement is central (i.e., that the
hyperplanes have a nonempty intersection); in this case $\F$ has an
identity.

The set $\C$ of chambers is an ideal, so we can run a random walk on
it.  A step in the walk goes from a chamber~$C$ to the chamber $FC$,
where $F$ is chosen according to some probability
distribution~$\{w_F\}_{F\in\F}$.  Examples in
\cite{\bid,\bhr,\bbd,\bd,\dicm} show that these \emph{hyperplane
chamber walks} include a wide variety of interesting processes.  The
references also explain a geometric interpretation of the step from
$C$ to~$FC$:  Namely, $FC$ is the chamber closest to~$C$ having $F$ as
a face.

Surprisingly, the eigenvalues of the transition matrix turn out to be
real.  In fact, they are certain partial sums of the weights $w_F$.
To say which partial sums occur, we need the \emph{intersection
lattice}~$\L$, consisting of all subspaces $X\subseteq V$ that are
intersections of some of the given hyperplanes; we order $\L$ by
inclusion.  The result, then, is that there is an eigenvalue
$\lambda_X = \sum_{F\subseteq X} w_F$ for each $X\in\L$, with
multiplicity $m_X=\abs{\mu(X,V)}$, where $\mu$ is the M\"{o}bius
function of~$\L$.  This was proved by Bidigare, Hanlon, and
Rockmore~\cite{\bhr}.  A different proof is given by Brown and
Diaconis \cite{\bd}, who show further that the transition matrix is
diagonalizable.

\subsection{Example:  The free LRB} \label{sub:free}

The \emph{free} LRB with identity on $n$ generators, denoted $F_n$,
may be constructed as follows: The elements of $F_n$ are sequences
$x=(x_1,\dots,x_l)$ of distinct elements of the set
$[n]=\{1,\dots,n\}$, $0\le l\le n$.  We multiply two such sequences by
\[
(x_1,\dots,x_l)(y_1,\dots,y_m)=(x_1,\dots,x_l,y_1,\dots,y_m)\sphat\,,
\]
where the hat means ``delete any element that has occurred earlier''.
For example,
\[
(2\,1)(3\,5\,4\,1\,6)=(2\,1\,3\,5\,4\,6).
\]
One can think of the elements of $F_n$ as reduced words on an alphabet
of $n$ letters, where ``reduced'' means that the word cannot be
shortened by applying~(D).  

The ideal $I$ on which we will run our random walk is the set of
reduced words of length~$n$; these can be identified with
permutations.  If the weights $w_x$ are concentrated on the $n$
generators, then the resulting random walk can be pictured as follows:
Think of $(x_1,\dots,x_n)$ as the set of labels on a deck of $n$
cards.  Then a step in the walk consists of removing the card labeled
$i$ with probability $w_i$ and replacing it on top.  This is the
well-studied \emph{Tsetlin library}, or \emph{weighted random-to-top
shuffle}, which arises in the study of dynamic list-management in
computer science.  See Fill~\cite{\Fill} and the references cited
there.

The eigenvalues were first found by Phatarfod~\cite{\phatarfod}; see
also \cite{\bid,\bhr,\bd,\Fill} for other proofs.

The result is that there is one eigenvalue $\lambda_X =\sum_{i\in X}
w_i$ for each subset $X\subseteq[n]$, with multiplicity equal to the
derangement number $d_{n-\abs{X}}$.  Here $d_k$ is the number of
fixed-point-free permutations of $k$ elements.  (Note that $d_1=0$, so
$\lambda_X$ does not actually occur as an eigenvalue if
$\abs{X}=n-1$.)

\subsection{Example:  A $q$-analogue} \label{sub:q-free}

Let $V$ be the $n$-dimensional vector space $\fq^n$, where $\fq$ is
the field with $q$-elements.  Let $F_{n,q}$ be the set of ordered
linearly independent sets $(x_1,\dots,x_l)$ in~$V$; two such are
multiplied by
\[
(x_1,\dots,x_l)(y_1,\dots,y_m)=(x_1,\dots,x_l,y_1,\dots,y_m)\sphat\,,
\]
where the hat means ``delete any vector that is linearly dependent on
the earlier vectors''.  Alternatively, we can think of the elements
of~$F_{n,q}$ as $n$-rowed matrices over~$\fq$ with independent
columns; we multiply two such matrices by juxtaposing them and then
deleting the columns that are linearly dependent on earlier columns.

A natural ideal to use is the set of ordered bases of $V$ or,
equivalently, the set of invertible matrices.  If we now assign
weights $w_v$ summing to~1 to the nonzero vectors $v\in V$ (i.e., to
the sequences $x$ as above of length~1), we get a Markov chain on
invertible matrices that can be described as follows:  Given an
invertible matrix, pick a nonzero vector $v$ with probability $w_v$
and adjoin it as a new first column; delete the unique column that is
linearly dependent on the earlier ones.

This chain does not seem to have been considered before.  We will see,
as a consequence of our main theorem, that its transition matrix is
diagonalizable, with an eigenvalue
\[
\lambda_X = \sum_{v\in X} w_v
\]
for each subspace $X\subseteq V$.  The multiplicity $m_X$ of this
eigenvalue is the number of elements of $\gl_n(\fq)$ with $X$ as
fixed subspace, i.e., the number of elements that fix $X$ pointwise
and act as a derangement on the set-theoretic complement $V - X$.

This Markov chain is, in some sense, a $q$-analogue of the Tsetlin
library.  We will construct in Section~\ref{s:q-tsetlin} a quotient
$\bar{F}_{n,q}$ of~$F_{n,q}$, for which the random walk is more
deserving of the name ``$q$-analogue of the Tsetlin library''.

\subsection{The main result} \label{sub:eigen}

If $S$ is any finite LRB with identity, one can construct an
associated lattice~$L$, along with a ``support map'' $\supp\colon
S\onto L$.  For the hyperplane face semigroup, $L$ is the intersection
lattice, the support of a face being its affine span.  For $S=F_n$,
$L$ is the lattice of subsets of~$[n]$, and the support of a word
$(x_1,\dots,x_l)$ is the underlying set $\{x_1,\dots,x_l\}$ of
letters.  And for $S=F_{n,q}$, $L$ is the lattice of subspaces
of~$\fq^n$, the support of $(x_1,\dots,x_l)$ being the subspace
spanned by $\{x_1,\dots,x_l\}$.  The ideal on which we run our random
walk is the set $C$ of all $c\in S$ with $\supp c=\onehat$, where
$\onehat$ is the largest element of~$L$.  Borrowing terminology from
the hyperplane example, we call the elements of~$C$ \emph{chambers}.
Our main result, illustrated by the examples above, can be stated
roughly as follows:

\begin{main}[Informal statement]
The transition matrix of the walk on chambers is diagonalizable, with
one eigenvalue
\begin{equation*} 
\lambda_X = \sum_{\supp y\le X} w_y
\end{equation*}
for each $X\in L$.  The multiplicity $m_X$ of this eigenvalue depends
on combinatorial data derived from $S$ and~$L$.
\end{main}

Unfortunately, the formula for $m_X$ is somewhat technical.  See
Theorem~\ref{t:main} in Section~\ref{s:main} for the precise
statement.

\subsection{Stationary distribution and convergence rate} \label{sub:stationary}

For the hyperplane chamber walk, Brown and Diaconis~\cite{\bd}
describe the stationary distribution and estimate the rate of
convergence to stationarity.  These results and their proofs apply
without change to the present setup.  For completeness, we state the
results here.  Note first that we can run, along with our walk on the
chambers, a random walk on $S$ starting at the identity; after $m$
steps it is at $x_m\cdots x_2 x_1$, where $x_1,x_2,\dots$ are i.i.d.\
with distribution $\{w_x\}$.  If $S$ is generated by $\{x\in
S:w_x\ne0\}$, then this walk is eventually in~$C$ with
probability~1.  Let $T$ be the first time $m$ that $x_m\cdots x_2
x_1\in C$.

\setcounter{theorem}{-1}

\begin{theorem} \label{t:stationary}
Let $S$ be a finite LRB with identity, and let $L$ be the associated
lattice.  Let $\{w_x\}$ be a probability distribution on~$S$ such that
$S$ is generated by $\{x\in S:w_x\ne0\}$.  Then the random walk on the
ideal $C$ of chambers has a unique stationary distribution~$\pi$; it
is the distribution of the infinite product $c=x_1x_2\cdots$, where
$x_1,x_2,\dots$ are i.i.d.\ with distribution~$\{w_x\}$.  The total
variation distance from stationarity after $m$ steps for the walk
started at any chamber $c_0$ satisfies
\[
\norm{P_{c_0}^m -\pi}\le \Pr\{T>m\} \le \sum_H \lambda_H^m,
\]
where $H$ ranges over the maximal elements (co-atoms) of $L$.
\end{theorem}

The fact that the infinite product converges (i.e., that the partial
sums are eventually constant) is an immediate consequence of
property~(D).  See \cite[Section~3]{\bd} for other descriptions of
$\pi$, involving sampling without replacement, that can be obtained by
making systematic use of~(D).

\subsection{Organization of the paper}

We begin by restarting the theory of LRBs in Section~\ref{s:lrb},
adopting a definition slightly different from (but equivalent to) the
one in Section~\ref{sub:semi}.  This allows us to get more quickly to
the main ideas of the paper without getting bogged down in semigroup
theory.  We can then give in Section~\ref{s:main} the precise
statement of our main theorem, with the multiplicities~$m_X$ spelled
out.  We also give some easy examples in that section.  

Sections \ref{s:convex}, \ref{s:q-tsetlin}, and \ref{s:matroid}
contain more elaborate examples.  Readers who wish to proceed to the
proof of the main theorem may skip ahead to Section~\ref{s:eigen}.  In
Section~\ref{s:convex} we consider a convex, open, polyhedral subset
$U\subset\R^n$, divided into chambers by hyperplanes that cut across
it.  There is a random walk on these chambers, generalizing the walk
of Section~\ref{sub:bhr}.  From a technical point of view, this is a
fairly trivial generalization; but it leads to new examples, including
a random walk on the maximal chains of any finite distributive
lattice.  An amusing special case of this is the ``kids walk''.
Section~\ref{s:q-tsetlin} treats the $q$-analogue of the Tsetlin
library mentioned above.  The multiplicities $m_X$ for this walk are
the $q$-derangement numbers studied by Wachs~\cite{\wachs}.  And
Section~\ref{s:matroid} gives a matroid generalization of the Tsetlin
library, including both the Tsetlin library and its $q$-analogue.
Applying the theory to graphical matroids, we obtain a random walk on
the edge-ordered spanning trees of a graph, as well as a random walk
that has a (speculative) connection with phylogenetic trees.

In Section~\ref{s:eigen} we begin the proof of the main theorem.  We
find the radical and semisimple quotient of the semigroup algebra $\R
S$ using ideas of Bidigare~\cite{\bid}, and from this we can read off
the irreducible representations of~$S$.  The eigenvalue formulas
follow easily, but not the diagonalizability of the transition matrix.

Diagonalizability is deduced in Section~\ref{s:semisimple} from a more
precise result, asserting that the subalgebra $\R[w]\subseteq\R S$
generated by~$w=\sum_{x\in S}w_x x$ is split semisimple (isomorphic to
a direct product of copies of~$\R$).  We use here a criterion that
deserves to be better known, involving the poles of the generating
function for the powers of~$w$.  As a byproduct of the proof, we
obtain an explicit (though complicated) formula for the primitive
idempotents in~$\R[w]$, and hence for the projections onto the
eigenspaces of~$P$.  Our methods are inspired by the work of
Fill~\cite{\Fill} on the Tsetlin library, and some of our formulas may
be essentially the same as unpublished results of his.  Finally, we
specialize in Section~\ref{s:descent} to the hyperplane face semigroup
of a reflection arrangement, and we give connections with Solomon's
descent algebra.  Here again we make crucial use of results of
Bidigare~\cite{\bid}.

There are three appendices that provide supplementary material.
Appendix \ref{app:hyperplane} summarizes the facts about hyperplane
arrangements that we use.  This appendix is not logically necessary,
but it is cited in many examples and it provides the motivation for
several definitions that would otherwise seem quite strange.
Appendix~\ref{app:semi} lays the foundations for the theory of LRBs;
in particular, it is here that we reconcile the definition given in
Section~\ref{sub:semi} with the one in Section~\ref{s:lrb}.  Finally,
in Appendix~\ref{app:derangement} we discuss a generalization of the
derangement numbers.  These arise naturally in connection with the
matroid examples of Section~\ref{s:matroid}.

\subsection*{Acknowledgments}

This paper grew out of my joint work with Persi Diaconis~\cite{\bd}.
It is a pleasure to thank him for many conversations and probing
questions, from which I have benefited enormously.  I am also grateful
to Swapneel Mahajan, Richard Stanley, and Michelle Wachs for helping
me with Appendix~\ref{app:derangement}.  John Howie provided helpful
pointers to the semigroup literature.  Finally, I would like to
acknowledge that, as is already evident, this paper owes a great debt
to the work of Bidigare~\cite{\bid}.

\subsection*{Convention}

For simplicity, all semigroups are assumed to be finite and to
have an identity, unless the contrary is stated explicitly.

\section{Left-regular bands} \label{s:lrb}

Let $S$ be a semigroup (finite, with identity).  It is convenient
to redefine ``LRB'' so that the lattice $L$, whose existence was
asserted in Section~\ref{sub:eigen}, is built into the definition.  The
interested reader can refer to Appendix~\ref{app:semi} for a proof
that the present definition is equivalent to the one in
Section~\ref{sub:semi}, as well as for other characterizations of LRBs.

\subsection{Definition} \label{sub:defn}

We say that $S$ is a \emph{LRB} if there are a lattice $L$ and a
surjection $\supp\colon S\onto L$ satisfying
\begin{equation} \label{e:homom}
\supp xy = \supp x \vee \supp y
\end{equation}
and
\begin{equation} \label{e:delete}
xy=x \quad\text{if } \supp y\le \supp x.
\end{equation}
Here $\vee$ denotes the join operation (least upper bound) in $L$.  It
follows from these axioms that every $x\in S$ is idempotent ($x^2=x$)
and that $S$ satisfies the identity
\begin{equation} \label{e:regular}
xyx=xy
\end{equation}
for all $x,y\in S$.  Thus $S$ has the deletion property~(D) stated in
Section~\ref{sub:semi}.

The motivation for \eqref{e:homom} and \eqref{e:delete} comes from the
theory of hyperplane arrangements (Appendix~\ref{app:hyperplane});
this theory, then, provides a huge supply of examples, one of which is
discussed in detail in Section~\ref{sub:braidsemi}.  Further examples
have been given in Sections \ref{sub:free} and~\ref{sub:q-free}, and
many more will be given in Sections \ref{s:convex}, \ref{s:q-tsetlin},
and~\ref{s:matroid},

\subsection{Partial order} \label{sub:poset}

If $S$ is a LRB, we can define a partial order on~$S$ by
setting
\begin{equation} \label{e:poset1}
x\le y \iff xy=y.
\end{equation}
(For motivation, see equation~\eqref{e:poset} in
Appendix~\ref{app:hyperplane}.)  This relation is reflexive because
every element of $S$ is idempotent.  And it is transitive because if
$xy=y$ and $yz=z$, then $xz=x(yz)=(xy)z=yz=z$.  To check antisymmetry,
suppose $x\le y$ and $y\le x$.  Then $xy=y$ and $yx=x$, hence
$x=yx=(xy)x=xy=y$, where the second-to-last equality
uses~\eqref{e:regular}; so $S$ is indeed a poset.

Note that left multiplication by $x$ is a \emph{projection}
(idempotent operator) mapping $S$ onto $S_{\ge x} = \{y\in S : y\ge
x\}$.  The latter is a LRB in its own right, the
associated lattice being the interval $[X,\onehat]$ in~$L$, where
$X=\supp x$ and $\onehat$ is the largest element of~$L$.  Note also
that $S_{\ge x}$ depends only on $X$, up to isomorphism.  Indeed, if
we also have $\supp x'=X$, then the projections (left multiplications)
defined by $x$ and $x'$ give mutually inverse semigroup isomorphisms
between $S_{\ge x}$ and~$S_{\ge x'}$; this is a straightforward
consequence of the axioms.  We may therefore write $S_{\ge X}$ instead
of $S_{\ge x}$.  Thus
\[
S_{\ge X} \iso \{y\in S : y\ge x\}
\]
for any fixed $x$ with $\supp x=X$.  Note that the random walk studied
in this paper is defined in terms of the projection operators
restricted to chambers, mapping $C$ onto $C_{\ge x}$.  For the
hyperplane face semigroup these projections have a geometric meaning
that we mentioned in Section~\ref{sub:bhr}.

Finally, we remark that there is also a LRB
\[
S_{\le X}=\{y\in S : \supp y\le X\},
\]
whose associated lattice is the interval $[\zerohat,X]$, where
$\zerohat$ is the smallest element of~$L$.

\subsection{Example:  The semigroup of ordered partitions} \label{sub:braidsemi}

One of the standard examples of a hyperplane arrangement is the braid
arrangement, which is discussed in detail in
\cite{\bid,\bhr,\bbd,\bd}; see also Section~\ref{sub:braid} of the
present paper.  Its face semigroup~$\B$ is easy to describe
combinatorially, without reference to hyperplane arrangements:  The
elements of~$\B$ are ordered partitions $B=(B_1,\dots,B_l)$ of the set
$[n]=\{1,2,\dots,n\}$.  Thus the $B_i$ are nonempty sets that
partition $[n]$, and their order matters.  We multiply two ordered
partitions by taking intersections and ordering them
lexicographically; more precisely, if $B=(B_1,\dots,B_l)$ and
$C=(C_1,\dots,C_m)$, then
\[
BC=(B_1\cap C_1,\dots,B_1\cap C_m,\dots, B_l\cap C_1,\dots,B_l\cap
C_m)\sphat\,,
\]
where the hat means ``delete empty intersections''.  This product
makes $\B$ a LRB, with the 1-block ordered partition as identity.  The
associated lattice $\L$ is the lattice of unordered set partitions
$\Pi$, with $\Pi\le\Pi'$ if $\Pi'$ is a refinement of~$\Pi$.  Thus the
smallest element $\zerohat$ of~$\L$ is the 1-block partition, and the
largest element $\onehat$ is the partition into singletons.  The
support map $\B\onto\L$ forgets the ordering of the blocks.

The partial order on $\B$ is also given by refinement, taking account
of the block ordering.  Thus $B\le C$ if and only if $C$ consists of
an ordered partition of $B_1$ followed by an ordered partition
of~$B_2$, and so on.  The chambers are the ordered partitions into
singletons, so they correspond to the permutations of~$[n]$.

It is useful to have a second description of $\B$.  Ordered partitions
$(B_1,\dots,B_l)$ of~$[n]$ are in 1--1 correspondence with chains of
subsets $\emptyset=E_0<E_1<\cdots<E_l=[n]$, the correspondence being
given by $B_i=E_i-E_{i-1}$.  So we may identify $\B$ with the set of
such chains.  The product is then described as follows:  Given a chain
$E$ as above and a second chain $F\colon\emptyset=F_0<F_1<\cdots<
F_m=[n]$, their product $EF$ is obtained by using $F$ to refine~$E$.
More precisely, consider the sets $G_{ij}=(E_{i-1}\cup F_j)\cap E_i =
E_{i-1}\cup (F_j\cap E_i)$.  For each $i=1,2,\dots,l$ we have
\[
E_{i-1} = G_{i0}\subseteq G_{i1}\subseteq\cdots\subseteq G_{im}=E_i.
\]
Deleting repetitions gives a chain from $E_{i-1}$ to~$E_i$, and
combining these for all~$i$ gives the desired refinement $EF$ of~$E$.

This construction is used in one of the standard proofs of the
Jordan--H\"{o}lder theorem.

\section{Statement of the main theorem} \label{s:main}

We are now in a position to complete the statement of our main result
by spelling out the multiplicities $m_X$ mentioned in
Section~\ref{sub:eigen}.  Let $S$ be a LRB with lattice of
supports~$L$.  For each $X\in L$ let $c_X$ be the number of chambers
in~$S_{\ge X}$, i.e., the number of chambers $c\in C$ such that $c\ge
x$, where $x$ is any fixed element of~$S$ having support~$X$.  Our
main theorem is:

\begin{theorem} \label{t:main}
Let $S$ be a finite LRB with identity, let $\{w_x\}$ be a probability
distribution on~$S$, and let $P$ be the transition matrix of the
random walk on chambers:
\begin{equation} \label{e:transition}
P(c,d) = \sum_{xc=d} w_x
\end{equation}
for $c,d\in C$.  Then $P$ is diagonalizable.  It has an eigenvalue
\[
\lambda_X = \sum_{\supp y\le X} w_y
\]
for each $X\in L$, with multiplicity $m_X$, where
\begin{equation} \label{e:mult}
\sum_{Y\ge X} m_Y = c_X
\end{equation}
for each $X\in L$.  Equivalently,
\begin{equation} \label{e:inv}
m_X = \sum_{Y\ge X} \mu(X,Y) c_Y,
\end{equation}
where $\mu$ is the M\"{o}bius function of the lattice $L$.
\end{theorem}

Note that $m_X$ depends only on the semigroup $S_{\ge X}$.  With this
in mind, there is an easy way to remember the formula~\eqref{e:mult}.
It says that for the random walk generated by any set of weights on
$S_{\ge X}$, the sum of the multiplicities of the eigenvalues is equal
to the number of chambers.

Here are a few easy examples.  More complicated examples will be
discussed in Sections \ref{s:convex}, \ref{s:q-tsetlin},
and~\ref{s:matroid}.

\begin{example} \label{ex:bhr}
Consider the chamber walk associated with a central hyperplane
arrangement in a vector space~$V$.  We have $c_X = \sum_{Y\ge X}
\abs{\mu(Y,V)}$ by Zaslavsky~\cite{\zas}.  Comparing this with
\eqref{e:mult}, we conclude that $m_X=\abs{\mu(X,V)}$.  Thus
Theorem~\ref{t:main} gives the results cited in Section~\ref{sub:bhr}.
The same results remain valid for noncentral arrangements.  This was
already shown in~\cite{\bd} by different methods.  To see how it
follows from Theorem~\ref{t:main}, one argues exactly as in the
central case, with one complication:  The face semigroup $\F$ need not
have an identity, and the poset $\L$ of supports of the faces is only
a semilattice (it has least upper bounds but not necessarily greatest
lower bounds).  Before applying the theorem, one has to adjoin an
identity to~$\F$ to get a LRB $\hat{\F}$ with identity, and one has to
adjoin a smallest element~$\zerohat$ to~$\L$ to get a
lattice~$\hat{\L}$.  The theorem would seem, then, to give an extra
eigenvalue~$\lambda_{\zerohat} =0$.  But Zaslavsky~\cite{\zas} showed
that $c_{\zerohat}$, the total number of chambers, is $\sum_{Y\in\L}
\abs{\mu(Y,V)}$.  One can now deduce from \eqref{e:mult} that
$m_{\zerohat}=0$ and hence that $\lambda_{\zerohat}$ does not really
occur as an eigenvalue.
\end{example}

\begin{example} \label{ex:free}
Let $S=F_n$.  As we stated in Section~\ref{sub:free}, $m_X$ is the
derangement number $d_{n-\abs{X}}$ for any $X\subseteq[n]$.  To deduce
this from Theorem~\ref{t:main}, we need only observe that
\begin{equation} \label{e:mult17}
\sum_{Y\supseteq X} d_{n-\abs{Y}} = c_X
\end{equation}
for each $X\subseteq[n]$.  Indeed, one can check that
$c_X=(n-\abs{X})!$, which is the number of permutations of~$[n]$ that
fix~$X$ pointwise; and the left-hand-side of~\eqref{e:mult17} counts
these according to their fixed-point sets.
\end{example}

\begin{example}
Let $S= F_{n,q}$.  We claimed in Section~\ref{sub:q-free} that $m_X$
for a subspace $X\subseteq\fq^n$ is the number of elements of
$\gl_n(\fq)$ with $X$ as fixed subspace.  To see this, note that $c_X$
is the number of ways to extend a given ordered basis of $X$ to an
ordered basis of $\fq^n$.  This is also the number of elements of
$\gl_n(\fq)$ that fix $X$ pointwise, and the claim now follows
from~\eqref{e:mult} exactly as in Example~\ref{ex:free}.
\end{example}

\section{Examples: Convex sets, distributive lattices, and the kids walk}
\label{s:convex}

The examples in this section were first treated in unpublished joint
work with Persi Diaconis, using the techniques of~\cite{\bd} rather
than semigroup methods.

Let $U\subset\R^n$ be a nonempty set that is a finite intersection of
open halfspaces.  A finite set of hyperplanes cutting across $U$
divides $U$ into regions.  We are interested in a random walk on these
regions driven by a set of weights on their faces.  A convenient way
to set this up is to combine the hyperplanes defining $U$ with the
hyperplanes cutting across $U$; this yields a hyperplane
arrangement~$\A$, and the regions into which $U$ is cut form a subset
$\D$ of the set~$\C$ of chambers of~$\A$.  Section~\ref{sub:convex}
spells out this point of view in more detail.  We then construct and
analyze a random walk on~$\D$ in Section~\ref{sub:walk}.  We show in
Section~\ref{sub:distrib} how the theory yields a random walk on the
maximal chains of a distributive lattice, and we illustrate this in
Section~\ref{sub:kids} by discussing the ``kids walk''.

\subsection{Convex sets of chambers} \label{sub:convex}

Let $\A=\{H_i\}_{i\in I}$ be a hyperplane arrangement in a real vector
space~$V$, let $\F$ be its face semigroup, and let $\C$ be the ideal
of chambers.  We do \emph{not} assume that $\A$ is central, so $\F$
need not have an identity.  Let $\D\subseteq\C$ be a convex set of
chambers, as defined in Section~\ref{sub:gallery}.  Thus there is a
subset $J\subseteq I$ and a set of signs $\sigma_i\in\{+,-\}$ ($i\in
J$) such that
\[
\D=\{C\in\C : \sigma_i(C)=\sigma_i \text{ for all } i\in J\}.
\]
We may assume that each $\sigma_i=+$.  The open set $U$ referred to
above is then $\bigcap_{i\in J} H_i^+$.

\begin{figure}[htb]
\begin{center}
\resizebox{!}{2in}{\includegraphics{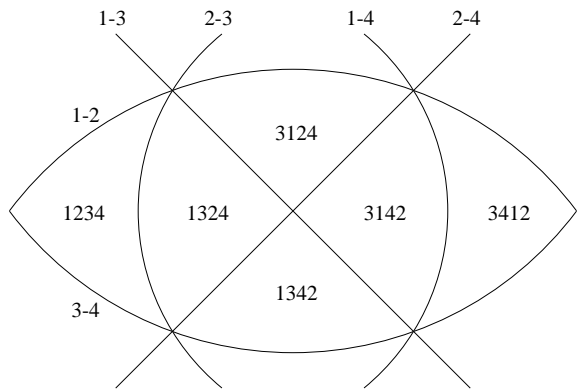}}
\end{center}
\caption{A convex subset of the braid arrangement.} \label{f:convex}
\end{figure}

As a simple example, consider the braid arrangement in $\R^4$
(Section~\ref{sub:braid}).  The region $U$ defined by $x_1>x_2$ and
$x_3>x_4$ contains six chambers, corresponding to the permutations
1234, 1324, 1342, 3124, 3142, 3412.  As explained in
Section~\ref{sub:cell}, it is possible to represent the arrangement by
means of a picture on the 2-sphere.  In this picture
(Figure~\ref{f:braid} in Section~\ref{sub:cell}) $U$ corresponds to
one of the open lunes determined by the great circles 1-2 and 3-4.
Figure~\ref{f:convex} gives a better view of this lune.

\subsection{A walk on the chambers} \label{sub:walk}

Let $\G$ be the set of faces of the chambers $D\in\D$; equivalently,
\[
\G=\{G\in\F : \sigma_i(G)\ge0 \text{ for all } i\in J\}.
\]
(To see that the right side is contained in the left, suppose
$\sigma_i(G)\ge0$ for all $i\in J$.  Choose an arbitrary $D\in\D$.
Then we have $G\le GD \in\D$, hence $G\in\G$.)  Then $\G$ is a
subsemigroup of~$\F$, hence a LRB (possibly without identity) in its
own right.  Its set of chambers is~$\D$.  Thus we can run a random
walk on~$\D$ driven by a set of weights on~$\G$.  To describe the
eigenvalues, we need some further notation.

Let $\G_0=\{G\in\G : \sigma_i(G)=+\}$.  In other words, $\G_0$ is the
set of faces that are contained in our open set~$U=\bigcap_{i\in J}
H_i^+$.  Let $\L$ be the intersection semilattice of~$\A$, let
$\M\subseteq\L$ be the set of supports of the faces in~$G$, and let
$\M_0\subseteq\M$ be the set of supports of the faces in~$\G_0$.
Equivalently, $\M_0$ consists of the $X\in\L$ that intersect~$U$.  In
our braid arrangement example, where we identify $\F$ with the set of
cells in the spherical representation of the arrangement, $\G_0$
consists of the cells in the interior of the lune:  one vertex, six
edges, and six chambers.  The bigger semigroup~$\G$ contains, in
addition, the six vertices and six edges on the boundary of the lune,
as well as the empty cell (which is the identity of~$\G$).  The
poset~$\M_0$ is shown in Figure~\ref{f:poset}.
\begin{figure}[htb]
\begin{center}
\setlength{\unitlength}{3947sp}%
\begingroup\makeatletter\ifx\SetFigFont\undefined
\def\x#1#2#3#4#5#6#7\relax{\def\x{#1#2#3#4#5#6}}%
\expandafter\x\fmtname xxxxxx\relax \def\y{splain}%
\ifx\x\y   
\gdef\SetFigFont#1#2#3{%
  \ifnum #1<17\tiny\else \ifnum #1<20\small\else
  \ifnum #1<24\normalsize\else \ifnum #1<29\large\else
  \ifnum #1<34\Large\else \ifnum #1<41\LARGE\else
     \huge\fi\fi\fi\fi\fi\fi
  \csname #3\endcsname}%
\else
\gdef\SetFigFont#1#2#3{\begingroup
  \count@#1\relax \ifnum 25<\count@\count@25\fi
  \def\x{\endgroup\@setsize\SetFigFont{#2pt}}%
  \expandafter\x
    \csname \romannumeral\the\count@ pt\expandafter\endcsname
    \csname @\romannumeral\the\count@ pt\endcsname
  \csname #3\endcsname}%
\fi
\fi\endgroup
\resizebox{!}{1.5in}{%
\begin{picture}(2487,2535)(376,-2386)
\thinlines
\put(451,-961){\line( 4, 3){1200}}
\put(1651,-61){\line( 4,-3){1200}}
\put(1651,-61){\line(-1,-2){450}}
\put(1651,-61){\line( 1,-2){450}}
\put(1201,-1261){\line( 1,-2){450}}
\put(2101,-1261){\line(-1,-2){450}}
\put(326,-1186){\makebox(0,0)[lb]{\smash{\SetFigFont{12}{14.4}{rm}$H_{14}$}}}
\put(1576, 14){\makebox(0,0)[lb]{\smash{\SetFigFont{12}{14.4}{rm}$V$}}}
\put(1100,-1186){\makebox(0,0)[lb]{\smash{\SetFigFont{12}{14.4}{rm}$H_{13}$}}}
\put(2026,-1186){\makebox(0,0)[lb]{\smash{\SetFigFont{12}{14.4}{rm}$H_{24}$}}}
\put(2776,-1186){\makebox(0,0)[lb]{\smash{\SetFigFont{12}{14.4}{rm}$H_{23}$}}}
\put(1226,-2386){\makebox(0,0)[lb]{\smash{\SetFigFont{12}{14.4}{rm}$H_{13}\cap H_{24}$}}}
\end{picture}}
\end{center}
\caption{The poset $\M_0$.} \label{f:poset}
\end{figure}

Note that if $X\in\M_0$ and $X\le Y\in\L$, then $Y\in\M_0$; this
implies that we get the same value for the M\"{o}bius number
$\mu(X,V)$ for $X\in\M_0$ no matter which of the posets $\M_0,\M,\L$
we work in.  We can now state:

\begin{theorem} \label{t:convex}
Let $\A$ be a hyperplane arrangement and let $\G$, $\D$, and~$\M_0$
be as above.  For any probability distribution $\{w_G\}_{G\in\G}$
on~$\G$, the transition matrix of the random walk on~$\D$ is
diagonalizable.  It has an eigenvalue
\[
\lambda_X = \sum_{\substack{G\in\G\\ G\subseteq X}} w_G
\]
for each $X\in \M_0$, with multiplicity $\abs{\mu(X,V)}$.
\end{theorem}

\begin{proof}
We argue as in our discussion of the walk on~$\C$ in
Example~\ref{ex:bhr} of Section~\ref{s:main}.  Assume first that $\A$
is central, so that $\G$ has an identity.  The lattice associated
with~$\G$ is~$\M$, so Theorem~\ref{t:main} gives us an
eigenvalue~$\lambda_X$ as above for each $X\in\M$, with
multiplicities~$m_X$ characterized by
\begin{equation} \label{e:mult8}
\sum_{\substack{Y\in\M\\Y\supseteq X}} m_Y = c_X
\end{equation}
for each $X\in \M$, where $c_X$ is the number of chambers in~$\G_X$.
We wish to show that $m_X=\abs{\mu(X,V)}$ for $X\in\M_0$ and that
$m_X=0$ for $X\notin\M_0$.  This will follow from \eqref{e:mult8} if
we show
\begin{equation} \label{e:zas}
\sum_{\substack{Y\in\M_0\\Y\supseteq X}} \abs{\mu(Y,V)} = c_X
\end{equation}
for each $X\in\M$.  

Now Zaslavsky \cite{\zasgen} counted the number of regions obtained when
an open convex set is cut by hyperplanes (see his Theorem~3.2 and the
comments at the bottom of p.~275).  His result, in our notation, is
\begin{equation} \label{e:zas1}
\abs{\D}=\sum_{Y\in\M_0} \abs{\mu(Y,V)}.
\end{equation}
This is the case $X=\zerohat$ of \eqref{e:zas}.  Equation~\eqref{e:zas}
for arbitrary~$X$ can be obtained by applying~\eqref{e:zas1} with $\A$
replaced by the set of hyperplanes $H\in\A$ that contain~$X$.  
Theorem~\ref{t:convex} is now proved if $\A$ is central.  

The noncentral case is treated by adjoining an identity to~$\G$, as in
Example~\ref{ex:bhr} of Section~\ref{s:main}.  The essential point is
that \eqref{e:zas1} still holds, and this implies that the ``extra''
eigenvalue $\lambda_{\zerohat} = 0$ has multiplicity~0.
\end{proof}

To illustrate the theorem, we return to the convex set in
Figure~\ref{f:convex}, with $\M_0$ as in Figure~\ref{f:poset}.  We
have $\mu(X,V)=\pm1$ for each $X\in\M_0$, so each contributes an
eigenvalue of multiplicity~1.  Suppose, for example, that we take
uniform weights $w_G=1/7$ on the seven vertices in
Figure~\ref{f:convex}.  Then Theorem~\ref{t:convex} gives the
following eigenvalues:
\[
\renewcommand{\arraystretch}{1.25}
\begin{array}{c|c}
X&\lambda_X\\ \hline
V&1\\
H_{13},H_{24}&3/7\\
H_{14},H_{23}&2/7\\
H_{13}\cap H_{24}&1/7
\end{array}
\]
The transition matrix $P$ in this case is $1/7$ of the following matrix:
\[
\begin{array}{c|cccccc}
&1234&1324&1342&3124&3142&3412\\ \hline
1234&3&1&1&1&0&1\\
1324&1&3&1&1&0&1\\
1342&1&1&3&0&1&1\\
3124&1&1&0&3&1&1\\
3142&1&0&1&1&3&1\\
3412&1&0&1&1&1&3
\end{array}
\]

\subsection{Distributive lattices} \label{sub:distrib}

If $L$ is a finite distributive lattice, there is a LRB $S$ whose
elements are chains $\zerohat=x_0<x_1<\cdots<x_l=\onehat$.  To
construct the product of two such chains, we use the second factor to
refine the first, exactly as in the discussion at the end of
Section~\ref{sub:braidsemi}, where we treated the Boolean lattice of
subsets of~$[n]$.  A simple way to verify that $S$ is indeed a LRB is
to appeal to the well-known fact that $L$ can be embedded as a
sublattice of a Boolean lattice.  Moreover,
Abels~\cite[Proposition~2.5]{\abels} has described a way of
constructing an embedding which makes the set of chambers in~$S$
(i.e., the maximal chains in~$L$) correspond to a convex set of
chambers in the braid arrangement.  His embedding depends on a choice
of a ``fundamental'' maximal chain, which then corresponds to the
identity permutation.  We can therefore use the results of
Section~\ref{sub:walk} to analyze a random walk on the maximal chains
of~$L$, driven by weights on arbitrary chains.

\begin{figure}[htb]
\begin{center}

\scalebox{.8}{\includegraphics{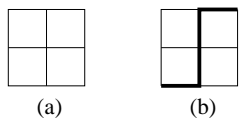}}
\end{center}
\caption{(a) A distributive lattice.  (b) A maximal chain.} \label{f:distrib}
\end{figure}

As an example of a distributive lattice, consider the product
$\{0,1,\dots,p\}\times \{0,1,\dots,q\}$ of a chain of length~$p$ by a
chain of length~$q$.  The case $p=q=2$ is shown in
Figure~\ref{f:distrib}(a).  The maximal chains are the lattice paths
from $(0,0)$ to~$(p,q)$, as in Figure~\ref{f:distrib}(b).  Each
maximal chain has length~$p+q$, and there are $\binom{p+q}{p}$ of
them; indeed, a lattice path can be identified with a binary vector of
length $p+q$ containing exactly $p$ ones.  (Think of 1 as ``right''
and 0 as ``up''.)

\begin{figure}[htb]
\begin{center}
\scalebox{.8}{\includegraphics{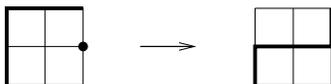}}
\end{center}
\caption{A step in the walk on lattice paths.} \label{f:pathmove}
\end{figure}

One interesting random walk on these lattice paths is obtained by
assigning uniform weights to the $(p+1)(q+1)-2$ chains of the form
$\zerohat<x<\onehat$.  A step in the walk consists of choosing $x\in
L-\{\zerohat,\onehat\}$ at random and then modifying the given path
minimally to make it pass through~$x$.  See Figure~\ref{f:pathmove} for
an illustration; here $x=(2,1)$.

\begin{figure}[htb]
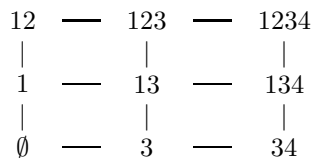

\[
\begin{matrix}
12&\lhor&123&\lhor&1234\\
|&&|&&|\\
1&\lhor&13&\lhor&134\\
|&&|&&|\\
\emptyset&\lhor&3&\lhor&34
\end{matrix}
\]
\caption{Embedding in the Boolean lattice.} \label{f:embed}
\end{figure}

In case $p=q=2$, the method of Abels cited above leads to an embedding
of $L=\{0,1,2\}\times\{0,1,2\}$ into the Boolean lattice of rank~4.
One such embedding is shown in Figure~\ref{f:embed}; it is obtained by
taking the fundamental maximal chain in~$L$ to be the lattice path
that goes up the left side and then across the top.  (Note:  An
expression like 134 in Figure~\ref{f:embed} denotes the set
$\{1,3,4\}$.)  The six maximal chains correspond to the six
permutations shown in Figure~\ref{f:convex}, and the walk on lattice
paths is the same as the walk discussed at the end of
Section~\ref{sub:walk}.

One can treat general $p,q$ in a similar way, but it would take us too
far afield to give further details.  One can also obtain results on
the stationary distribution and convergence rate via
Theorem~\ref{t:stationary} (Section~\ref{sub:stationary}).

\subsection{The kids walk} \label{sub:kids}

This is a walk on the $p$-subsets of a $(p+q)$-set, represented as
binary vectors of length~$p+q$ containing $p$ ones.  Think of the
zeroes as empty spaces and the ones as spaces occupied by kids.  At
each step a kid and an empty space are independently chosen at
random.  The kid then moves toward the empty space, pushing any other
kids he encounters until the space is occupied.  Here is an example
with $p=3$ and $q=4$.  The initial configuration is
\[
\begin{matrix}
&\downarrow&&&\downarrow\\
0&1&1&0&0&1&0
\end{matrix}
\]
with the two chosen positions indicated by arrows.  The final
configuration is
\[
\begin{matrix}
\\
0&0&0&1&1&1&0
\end{matrix}
\]

The kids walk is the same as the walk on lattice paths described in
Section~\ref{sub:distrib}, except that the latter has holding; namely,
the chosen lattice point~$x$ is on the current path with probability
$\alpha=(p+q-1)/(pq+p+q-1)$, in which case the walk stays at the
current path.  Thus if $P$ is the transition matrix for the walk on
lattice paths and $P_1$ is the transition matrix for the kids walk, we
have $P=\alpha I + (1-\alpha)P_1$, so that $P_1=(P-\alpha
I)/(1-\alpha)$.  It follows that $P_1$ is diagonalizable with
eigenvalues $(\lambda-\alpha)/(1-\alpha)$, where $\lambda$ ranges over
the eigenvalues of~$P$.  If $p=q=2$, for example, we have
$\alpha=3/7$, and the result at the end of Section~\ref{sub:walk}
gives eigenvalues $1,0,0,-1/4,-1/4,-1/2$ for the kids walk.

\section{Example: A $q$-analogue of the Tsetlin library} \label{s:q-tsetlin}

The random walk in this section is based on a quotient $\bar{F}_{n,q}$
of the semigroup~$F_{n,q}$ (Section~\ref{sub:q-free}).  For
motivation, we begin by defining a quotient $\bar{F}_n$ of~$F_n$, and
we explain how it is related to the Tsetlin library.  The $q$-analogue
is then given in Section~\ref{sub:q-analogue}.

\subsection{A quotient of $F_n$} \label{sub:fnbar}

The references cited in Section~\ref{sub:braidsemi} show how the
random walk associated with the semigroup~$\B$ of ordered partitions
captures many shuffling schemes.  In particular, to obtain the Tsetlin
library one puts weight $w_i>0$ on the 2-block ordered partition
$(i,[n]-i)$ and weight 0 on all other ordered partitions, where
$\sum_{i=1}^n w_i =1$.  From the point of view of the present paper,
however, the semigroup $\B$ is much too big for the study of the
Tsetlin library; one should replace $\B$ by the subsemigroup (with
identity) generated by the $n$ two-block ordered partitions to which
we have assigned weights.  This subsemigroup, which we denote by
$\bar{F}_n$, is easily described:  It consists of the ordered
partitions $(B_1,\dots,B_l)$ such that each block $B_i$ is a singleton
except possibly $B_l$.  Alternatively, it consists of the chains
$\emptyset=E_0<E_1<\cdots<E_l=[n]$ with $\abs{E_i}=i$ for $0\le i<l$.

The freeness of $F_n$ implies that $\bar{F}_n$ is a quotient of~$F_n$.
Explicitly we have a surjection $F_n\onto\bar{F}_n$ sending the
sequence $(x_1,\dots,x_l)$ to the following ordered partition~$B$:  If
$l<n$, then $B$ has $l+1$ blocks, with $B_i=\{x_i\}$ if $i\le l$ and
$B_{l+1}=[n]-\{x_1,\dots,x_l\}$; if $l=n$, then $B$ is the partition
into singletons $B_i=\{x_i\}$, $1\le i\le n$.  In terms of chains of
subsets, $B$ corresponds to the chain~$E$ with $E_i=\{x_1,\dots,x_i\}$
for $1\le i\le l$ and, if $l<n$, $E_{l+1}=[n]$.

The lattice of supports of~$\bar{F}_n$ can be identified with the set
of subsets $X\subseteq[n]$ such that $\abs{X}\ne n-1$, the support
of~$B$ being the union of the singleton blocks.  Note that the
quotient map $F_n\onto\bar{F}_n$ is almost 1--1; the only
identifications are that each $(n-1)$-tuple $(x_1,\dots,x_{n-1})$
in~$F_n$ gets identified with its (unique) extension to an $n$-tuple
in~$F_n$.

\begin{remark}
The semigroups $F_n$ and $\bar{F}_n$ have the same set of chambers,
and we have seen that either one can be used to generate the Tsetlin
library.  But $\bar{F}_n$ is more efficient for this purpose, in the
following two senses: (a)~When we use $F_n$, Theorem~\ref{t:main}
gives extraneous eigenvalues $\lambda_X$ with $\abs{X}=n-1$, which
then turn out not to occur because $m_X=0$.  (b)~The estimate of
convergence rate given in Theorem~\ref{t:stationary} is sharper if we
use $\bar{F}_n$ than if we use $F_n$, because the maximal elements of
the support lattice have size $n-2$ instead of~$n-1$.
\end{remark}

\subsection{$q$-analogue} \label{sub:q-analogue}

Let $V$ be the vector space $\fq^n$, where $\fq$ is the field with
$q$ elements.  As a $q$-analogue of~$\bar{F}_n$ we propose the
following semigroup~$\bar{F}_{n,q}$:  An element of~$\bar{F}_{n,q}$ is
a chain of subspaces $0=X_0<X_1<\cdots<X_l=V$ with $\dim X_i =i$ for
$i<l$.  Thus the chain cannot be refined except possibly at the last
step, between $X_{l-1}$ and~$V$.  Given two such chains
$\mathbf{X}=(X_0,\dots,X_l)$ and $\mathbf{Y}=(Y_0,\dots,Y_m)$, we
construct the product $\mathbf{X}\mathbf{Y}$ by using $\mathbf{Y}$ to
refine $\mathbf{X}$.  More precisely, the product is obtained by
forming the chain
\[
0=X_0<\cdots<X_{l-1}\le X_{l-1}+ Y_1 \le X_{l-1}+Y_2\le\cdots\le
X_{l-1}+ Y_m =V
\]
and deleting repetitions.  

The simplest way to verify that this product is associative is to
exhibit $\bar{F}_{n,q}$ as a quotient of the semigroup~$F_{n,q}$ of
ordered independent sets.  Namely, we can map $F_{n,q}$
onto~$\bar{F}_{n,q}$ by sending $(x_1,\dots,x_l)$ to the chain with
$X_i$ equal to the span of $\{x_1,\dots,x_i\}$ for $0\le i\le l$ and,
if $l<n$, $X_{l+1}=V$.  This gives a product-preserving surjection
$F_{n,q}\onto\bar{F}_{n,q}$, so our product on $\bar{F}_{n,q}$ is
indeed associative.

It is easy to check that $\bar{F}_{n,q}$ is a LRB whose associated
lattice is the set of subspaces of~$V$ of dimension different
from~$n-1$.  The support map is given by
\[
\supp (X_0,\dots,X_l)=
\begin{cases}
X_{l-1}& \text{if $l<n$}\\
V& \text{if $l=n$.}
\end{cases}
\]
Note that the join of two such subspaces $X,Y$ in this lattice is
their vector space sum $X+ Y$ unless the latter has dimension~$n-1$,
in which case the join is~$V$.

The chambers in $\bar{F}_{n,q}$ are the maximal chains
$0=X_0<X_1<\cdots<X_n=V$.  To construct a random walk analogous to the
Tsetlin library, put weight $w_\ell>0$ on the chain $0<\ell<V$ for
each subspace $\ell$ of dimension~1, and put weight 0 on all other
elements of~$\bar{F}_{n,q}$, where $\sum_\ell w_\ell=1$.  This yields
a walk on maximal chains that goes as follows:  Given a maximal chain
\[
0<X_1<\cdots<X_{n-1}<V,
\]
pick a line $\ell$  with probability $w_\ell$, and form a new
maximal chain
\[
0<\ell\le \ell+ X_1\le\cdots\le\ell+ X_{n-1}\le V;
\]
exactly one of the inequalities is an equality, and we delete the
repetition.  One can also view this walk as taking place on the
maximal flags in the projective space $\mathbb{P}^{n-1}(\fq)$, driven
by weights on the points.  If $n=3$, for example, this is a walk on
the incident point-line pairs in the projective plane.

According to Theorem~\ref{t:main} the transition matrix of this walk
is diagonalizable, with an eigenvalue
\[
\lambda_X = \sum_{\ell\subseteq X} w_\ell
\]
for each subspace $X$ with $\dim X\ne n-1$; the multiplicities $m_X$ are
characterized by
\begin{equation*}
\sum_{Y\supseteq X} m_Y = c_X,
\end{equation*}
where $c_X$ is the number of maximal chains in the interval
$[X,V]$ in~$L$.  It follows that $m_X$ is the $q$-derangement
number $d_{n-\dim X}(q)$ of Wachs~\cite{\wachs}; see
Example~\ref{ex:qder} in Section~\ref{sub:examples}.  This is why we
view the present walk as the ``right'' $q$-analogue of the Tsetlin
library, rather than the walk based on~$F_{n,q}$.

The stationary distribution $\pi$ of this walk is a probability
measure on the set of maximal chains.  One can deduce from
Theorem~\ref{t:stationary} the following description of~$\pi$:  Sample
from the set of lines~$\ell$ according to the weights $w_\ell$ to get
a line $l_1$.  Remove $\ell_1$ and sample again to get $\ell_2$.
Remove all the lines contained in $\ell_1+\ell_2$ and choose
$\ell_3$.  Continuing in this way, we obtain a maximal chain
\[
0<\ell_1<\ell_1+\ell_2<\cdots<\ell_1+\ell_2+\cdots+\ell_{n-1}<V
\]
after $n-1$ steps.  This chain is distributed according to~$\pi$.

\begin{remark}
This $q$-analogue of the Tsetlin library was first studied in joint
work with Persi Diaconis [unpublished], in which we extended the
hyperplane chamber walks to walks on the chambers of a building.  And
the first proof that the multiplicities were given by the
$q$-derangement numbers was arrived at with the help of Richard
Stanley.  In fact, the original calculation of the multiplicities,
which was quite different from the one given in this paper, led to
formulas similar to those of Proposition~\ref{p:defn}
(Section~\ref{sub:defn1}), and it was not immediately obvious that
these formulas gave the $q$-derangement numbers.
\end{remark}

\section{Examples: Random walks associated with matroids} \label{s:matroid}

Matroids were introduced by Whitney \cite{\whitney}, as an abstraction
of the linear independence properties of the columns of a matrix.  We
describe in this section two natural LRBs $S,\bar{S}$ that can be
associated with a matroid, hence two random walks.  These generalize
the pairs $F_n,\bar{F}_n$ and $F_{n,q},\bar{F}_{n,q}$ discussed in
Section~\ref{s:q-tsetlin}.

We begin by reviewing matroid concepts in Section~\ref{sub:review}.
We then construct the semigroups and the associated walks in
Section~\ref{sub:matroidsemi}.  Our discussion is brief because the
theory follows quite closely the two special cases already treated.
In Section~\ref{sub:graph} we consider a third case, graphical
matroids.  This leads to two random walks associated with a graph.  In
an effort to understand one of these examples intuitively, we give an
interpretation of it in terms of phylogenetic trees.

\subsection{Review of matroids} \label{sub:review}

The book by Welsh~\cite{\welsh} is a good reference for this
subsection.  A \emph{matroid} $M$ consists of a finite set $E$ and a
collection of subsets of~$E$, called \emph{independent sets}, subject
to axioms modeled on the notion of linear independence in vector
spaces.  A maximal independent set is called a \emph{basis}, and all
such have the same cardinality $n$, called the \emph{rank} of~$M$.
More generally one can define $\rank(A)$ for any subset $A\subseteq E$
as the rank of any maximal independent set in~$A$.  Any such
maximal independent set is called a basis for~$A$.  We say that
$x$ \emph{depends} on~$A$ if $\rank(A\cup x)=\rank(A)$; otherwise,
$\rank(A\cup x)=\rank(A)+1$.  Any set $A$ has a \emph{closure}
$\sigma(A)$, obtained by adjoining every $x$ that depends on~$A$, and
$A$ is said to be \emph{closed}, or a \emph{flat}, if $\sigma(A)=A$.
The set $L$ of flats is a lattice under inclusion.  One can think of
$L$ as an analogue of the lattice of subspaces of a vector space, and
$\sigma(A)$ plays the role of the span of a set of vectors.

In addition to the motivating example, in which $E$ is a set of
vectors, there are two other canonical examples:  The first is the
\emph{free} matroid of rank~$n$; the set $E$ is $\{1,2,\dots,n\}$, and
all subsets are independent.  The second is the \emph{graphical
matroid} associated with a finite graph~$G$; here $E$ is the set of
edges of~$G$, and a subset is independent if it contains no cycles.

\subsection{Semigroups associated with a matroid} \label{sub:matroidsemi}

Let $M$ be a matroid of rank $n$ with underlying set $E$.  We
construct two LRBs $S,\bar{S}$.  The elements of~$S$ are ordered
independent sets, i.e., $l$-tuples $\mathbf{x}=(x_1,\dots,x_l)$ of
distinct elements of~$E$ whose underlying set $\{x_1,\dots,x_l\}$ is
independent.  We set $\supp \mathbf{x}=\sigma(\{x_1,\dots,x_l\})$.
The product is defined by
\[
(x_1,\dots,x_l)(y_1,\dots,y_m)=(x_1,\dots,x_l,y_1,\dots,y_m)\sphat\,,
\]
where the hat means ``delete any element that depends on the earlier
elements''.  It is easy to check that we obtain in this way a LRB $S$
whose associated lattice is the lattice of flats~$L$.

The chambers of $S$ are the ordered bases of $M$.  To construct a
random walk analogous to the Tsetlin library, put weight $w_x>0$ on
the 1-tuple $(x)$ and weight 0 on the $l$-tuples with $l\ne1$, where
$\sum_x w_x =1$.  (Note: Not all $x\in E$ occur here, since $M$ might
contain \emph{loops}, i.e., elements $x$ such that the singleton
$\{x\}$ is not independent.)  This yields a walk on ordered bases that
goes as follows:  Given an ordered basis $(x_1,\dots,x_n)$, pick a
nonloop $x\in E$ with probability $w_x$, and make it the new first
basis element; delete the (unique) $x_i$ that depends on
$\{x,x_1,\dots,x_{i-1}\}$.

According to Theorem \ref{t:main}, the transition matrix of this walk
is diagonalizable, with an eigenvalue
\[
\lambda_X = \sum_{x\in X} w_x
\]
for each flat $X$; the multiplicities $m_X$ are characterized by
\begin{equation} \label{e:mult2}
\sum_{Y\ge X} m_Y = c_X
\end{equation}
for each $X\in L$.  Here $c_X$ is the number of ways of completing any
fixed basis of~$X$ to a basis of~$M$.  The stationary distribution
$\pi$ of this chain is a probability measure on the set of ordered
bases.  One can deduce from Theorem~\ref{t:stationary} the following
description of~$\pi$:  Sample from $E$ (according to the weights
$w_x$) to get a nonloop $x_1$.  Remove $x_1$ and everything dependent
on it and sample again to get $x_2$.  Remove the closure of
$\{x_1,x_2\}$, choose $x_3$, and so on.  After $n$ steps we  have
an ordered basis $(x_1,\dots,x_n)$ whose distribution is~$\pi$.

The second semigroup, $\bar{S}$, consists of chains of flats
$\zerohat=X_0<X_1<\cdots<X_l=\onehat$ with $\rank(X_i) =i$ for $i<l$.
Given two such chains $\mathbf{X}=(X_0,\dots,X_l)$ and
$\mathbf{Y}=(Y_0,\dots,Y_m)$, their product $\mathbf{X}\mathbf{Y}$ is
obtained by forming the chain
\[
\zerohat=X_0<\cdots<X_{l-1}\le X_{l-1}\vee Y_1 \le \cdots\le
X_{l-1}\vee Y_m =\onehat
\]
and deleting repetitions.  Here $\vee$ denotes the join operation in
the lattice of flats, i.e., $X\vee Y=\sigma(X\cup Y)$.  One can
verify, exactly as in Section~\ref{sub:q-analogue}, that $\bar{S}$ is
a LRB whose associated lattice $\bar{L}$ is the set of flats of rank
different from~$n-1$.

The chambers in $\bar{S}$ are the maximal chains
$\zerohat=X_0<X_1<\cdots<X_n=\onehat$.  To construct a random walk
analogous to the Tsetlin library, put weight $w_\ell>0$ on the chain
$\zerohat<\ell<\onehat$ for each flat $\ell$ of rank~1, where
$\sum_\ell w_\ell=1$.  This yields a walk on maximal chains that goes as
follows:  Given a maximal chain
\[
\zerohat<X_1<\cdots<X_{n-1}<\onehat,
\]
pick a flat $\ell$ of rank 1 with probability $w_\ell$, and form a new
maximal chain
\[
\zerohat<\ell\le \ell\vee X_1\le\cdots\le\ell\vee X_{n-1}\le\onehat;
\]
exactly one of the inequalities is an equality, and we delete the
repetition.

According to Theorem~\ref{t:main}, the transition matrix of this walk
is diagonalizable, with an eigenvalue
\[
\lambda_X = \sum_{\ell\le X} w_\ell
\]
for each flat $X$ with $\rank(X)\ne n-1$; the multiplicities $m_X$ are
characterized by
\begin{equation} \label{e:mult3}
\sum_{Y\ge X} m_Y = c_X
\end{equation}
where now $c_X$ is the number of maximal chains in the interval
$[X,\onehat]$ in~$L$.  

Recall that the multiplicities of the eigenvalues for the Tsetlin
library and its $q$-analogue are the derangement numbers and their
$q$-analogues.  Motivated by this, we show in
Appendix~\ref{app:derangement} how to associate a ``generalized
derangement number'' $d(L)$ to every finite lattice~$L$.  It will
follow quickly from the definition that the multiplicities
in~\eqref{e:mult3} are given by
\begin{equation} \label{e:multder}
m_X=d\left([X,\onehat]\right);
\end{equation}
see equation~\eqref{e:recurrencey} and the discussion following it.

The stationary distribution $\pi$ of this chain is a probability
measure on the set of maximal chains.  One can deduce from
Theorem~\ref{t:stationary} the following description of~$\pi$:  Sample
from the set of ``lines'' (rank~1 flats) according to the weights
$w_\ell$ to get a line $l_1$.  Remove $\ell_1$ and sample again to get
$\ell_2$.  Remove all the lines contained in $\ell_1\vee\ell_2$ and
choose $\ell_3$, and so on.  After $n-1$ steps we  have a maximal
chain
\[
\zerohat<\ell_1<\ell_1\vee\ell_2<\cdots<\ell_1\vee\ell_2\vee\cdots
\vee\ell_{n-1}<\onehat,
\]
which is distributed according to~$\pi$.

\subsection{Random walks associated with graphs} \label{sub:graph}

One of Whitney's main motivations in developing the theory of matroids
was the connection with graph theory.  As we mentioned in
Section~\ref{sub:review}, every finite graph $G$ gives rise to a
matroid whose underlying set $E$ is the set of edges of~$G$, with a
set of edges being independent if it contains no cycles.
Equivalently, the independent sets correspond to forests $F\subseteq
G$, where we make the convention that a forest contains every vertex
of~$G$.  We briefly describe here our two random walks, as specialized
to the matroid of~$G$.  Much remains to be understood about these
examples.

For simplicity, all of our graphs are assumed \emph{simple} (no
loops or multiple edges).  The lattice of flats $L=L(G)$ of the
graphical matroid can then be described as follows.  An element of $L$
is specified by a partition $\Pi$ of the vertex set $V$ of~$G$ such
that each block induces a connected subgraph.  The ordering on $L$ is
given by refinement, but with the opposite convention from the one
used in Section~\ref{sub:braidsemi}:  In $L(G)$, $\Pi\le\Pi'$ if
$\Pi$ is a refinement of~$\Pi'$.  Thus going up in the lattice
corresponds to merging blocks.  Associated with each $\Pi\in L$ is the
\emph{contraction} $\bar{G}=G/\Pi$, obtained by collapsing each block
to a point and making the resulting graph simple.  (Delete loops and
replace multiple edges by a single edge.)  Equivalently, $G/\Pi$ is
the simple graph with one vertex for each block, two blocks $B,B'$
being adjacent if, in~$G$, some vertex in~$B$ is adjacent to some
vertex in~$B'$.  Because of this interpretation of partitions, $L(G)$
is often called the \emph{lattice of contractions of $G$}.  The
smallest element $\zerohat$ is the partition into singletons (so
$\bar{G}=G$), and the largest element $\onehat$ is the partition into
connected components (so $\bar{G}$ is the discrete graph with one
vertex for each connected component of $G$).  From the collapsing
point of view, going up in the lattice $L(G)$ corresponds to doing
further collapsing.

Consider now the two semigroups $S,\bar{S}$ associated with the
graphical matroid.  An element of~$S$ can be identified with an
\emph{edge-ordered} forest $F\subseteq G$, i.e., a forest together
with a linear ordering of its edges.  The support of $F$ is the
partition of~$V$ given by the connected components of~$F$.  In
particular, the chambers of $S$ are the edge-ordered spanning forests
of~$G$ (spanning trees if $G$ is connected).  The random walk on these
chambers goes as follows:  Given a spanning forest with ordered edges
$e_1,\dots,e_n$, pick an edge $e$ with probability $w_e$ and make it
the new first edge; delete the first $e_i$ such that
$\{e,e_1,\dots,e_i\}$ contains a cycle.

We leave it to interested reader to spell out what the general results
in Section~\ref{sub:matroidsemi} say about this example.  One
interesting question arises:  Running this random walk, say with
uniform weights, gives a way of choosing an edge-ordered spanning
forest with distribution~$\pi$; what is the distribution of the
spanning forest obtained by forgetting the ordering?

We turn next to the random walk on maximal flags, based on the
semigroup~$\bar{S}$.  A maximal flag in $L(G)$ is gotten by collapsing
an edge of $G$ to get a (simple) graph~$G_1$, then collapsing an edge
of~$G_1$ to get~$G_2$, and so on, until we reach a discrete
graph~$G_n$.  The number $n$ of collapses is the number of edges in a
spanning forest of~$G$, i.e., the number of vertices of~$G$ minus the
number of connected components.  Note that an edge-ordered spanning
forest determines a collapsing sequence (maximal flag), but this
correspondence is not 1--1.  Different edge-ordered spanning forests
can give the same maximal flag, just as different ordered bases of a
vector space can determine the same maximal flag of subspaces.

We close this section by giving an interpretation of these
maximal flags and the corresponding random walk in terms of
phylogenetic trees.  Think of the vertices of~$G$ as species that
exist today.  We join two species by an edge if we think they might
have had a direct common ancestor.  Thus humans and chimpanzees are
probably adjacent, but not humans and frogs.  Assume, for simplicity,
that $G$ is connected, so that all species ultimately evolved from one
common ancestor.  To run the random walk, we are given weights on the
edges.  These can be thought of as indicating the strength of our
belief that two species have a direct common ancestor; alternatively,
they might indicate how recently we think they diverged from that
ancestor.

Recall that a maximal flag in~$L(G)$ consists of a sequence of edge
collapses
\[
G=G_0\onto G_1\onto\cdots\onto G_n=\text{point}.
\]
We can think of this as representing a feasible reconstruction of the
phylogenetic tree describing the evolution from the original common
ancestor to the present-day situation, in reverse chronological order.
Thus the first edge collapsed corresponds to the two species that most
recently split off from a direct common ancestor.  The collapsed graph
$G_1$ then represents the situation before that split.  The edge
of~$G_1$ that is collapsed to form~$G_2$ corresponds to the
next-most-recent split, and so on.

The random walk proceeds as follows:  Given a collapsing sequence as
above, pick a random edge $e$ of~$G$ according to the weights.  Make a
new collapsing sequence in which $e$ is collapsed first, but after
that the collapses mimic those of the original sequence.  In other
words, we revise our view of the evolutionary history by declaring
that two particular species were the most recent to split from a
common direct ancestor.

A pick from the stationary distribution of this walk can, as usual, be
obtained by sampling without replacement.  In the present situation
this amounts to the following: Pick an edge of~$G$ according to the
weights and collapse it to get $G_1$.  Use the collapsing map $G\onto
G_1$ to put weights on the edges of~$G_1$; thus the weight on an edge
of $G_1$ is the sum of weights of the edges of~$G$ that map to that
edge.  Note that the weights on $G_1$ do not sum to~1, because at
least one edge of~$G$ with positive weight gets collapsed to a point
in~$G_1$; so we must rescale them.  Now repeat the process:  Choose an
edge of $G_1$ according to the weights and collapse it to get $G_2$.
Continue in this way until a maximal collapsing sequence is obtained.

\begin{remark}
See Aldous~\cite{aldous99:_moran_markov} for a detailed analysis of
this walk in the case of uniform weights.
\end{remark}

\section{Irreducible representations and computation of eigenvalues} \label{s:eigen}

We now begin the proof of Theorem \ref{t:main}, starting with the
description of the eigenvalues.  Throughout this section $S$ denotes a LRB
(finite, with identity), and $\supp\colon S\onto L$ is the associated
support map.  Assume that we are given a probability distribution
$\{w_x\}$ on $S$ and that $P$ is the transition matrix of the random
walk on chambers.  We begin by recalling in Section~\ref{sub:algebra}
the algebraic interpretation of $P$ in terms of the semigroup algebra
of~$S$.  In Section~\ref{sub:radical} we compute the radical and
semisimple quotient of the semigroup algebra.  This was done by
Bidigare~\cite{\bid} for hyperplane face semigroups, and the proof in
general is identical.  We include the proof for the convenience of the
reader, since the thesis~\cite{\bid} is not readily available.  From
this result we can read off the irreducible representations of~$S$,
and the eigenvalue formula stated in Theorem~\ref{t:main} follows at
once; we explain this in Section~\ref{sub:eigen1}.

\subsection{Algebraic formulation} \label{sub:algebra}

It is well-known to probabilists that the transition matrix of a
random walk on a semigroup can be interpreted as the matrix of a
convolution operator.  (This is perhaps best known for groups, but the
result remains valid for semigroups.)  We wish to recast this result
in ring-theoretic language.  Consider the vector space $\R S$ of
formal linear combinations $\sum_{x\in S} a_x x$ of elements of~$S$,
with $a_x \in \R$.  The product on~$S$ extends to a bilinear product
on~$\R S$, making the latter a ring (the \emph{semigroup ring} of $S$
over~$\R$).  Thus
\[
\left(\sum_{x\in S} a_x x\right) \left(\sum_{x\in S} b_x x\right) =
\sum_{x\in S} c_x x,
\]
where
\[
c_x = \sum_{yz=x} a_y b_z.
\]
On the level of coefficients, this is the familiar \emph{convolution
product}.

A probability distribution $\{w_x\}_{x\in S}$ can be encoded in the
element
\[
w=\sum_{x\in S} w_x x
\]
of $\R S$, and I claim that the transition matrix $P$ of the random
walk determined by~$\{w_x\}$ is simply the matrix of the operator
``left multiplication by~$w$''.  More precisely, we have for any
$a=\sum_s a_s s$ in~$\R S$
\begin{align*}
wa &= \sum_x w_x x \sum_s a_s s = \sum_t \left(\sum_{\substack{x,s\\xs=t}}
w_x a_s\right) t\\
  &= \sum_t \left(\sum_s a_s P(s,t)\right)t,
\end{align*}
where the last equality follows from \eqref{e:transition1}.  Thus left
multiplication by $w$ acting on $\R S$ corresponds to right
multiplication by $P$ acting on row vectors $(a_s)_{s\in S}$.
Similarly, if we run the walk on an ideal $C\subseteq S$, then the
transition matrix is the matrix of left multiplication by $w$ on $\R
C$, which is an ideal in the ring~$\R S$.

In principle, then, the analysis of the random walk has been reduced
to ring theory.  Here is a familiar example in which this point of
view can be exploited (using $\mathbb{C}$ instead of~$\R$).  Suppose
that $S$ is a finite abelian group~$G$, and let $\hat{G}$ be its group
of characters $\chi\colon G\to\mathbb{C}^*$.  Then the Fourier
transform gives a ring isomorphism
\[
\mathbb{C}G \isoto \mathbb{C}^{\hat{G}}.
\]
Here $\mathbb{C}^{\hat{G}}$ is the ring of functions
$\hat{G}\to\mathbb{C}$ (with functions multiplied pointwise), and the
Fourier transform of $a=\sum_x a_x x$ is the function~$\hat{a}$ given
by $\hat{a}(\chi)=\sum_x a_x\chi(x)$; see~\cite[Section~6.2]{\serre}.
In particular, left multiplication by our element $w$ acting
on~$\mathbb{C}G$ is transformed to multiplication by~$\hat{w}$ acting
on~$\mathbb{C}^{\hat{G}}$.  This operator is diagonal with respect to
the standard basis of~$\mathbb{C}^{\hat{G}}$, and one concludes that
the eigenvalues of the transition matrix~$P$ are simply the
numbers~$\hat{w}(\chi)$.  Moreover, the Fourier inversion formula
gives an explicit diagonalization of multiplication by~$w$ and hence
of~$P$.

\subsection{Structure and representations of the semigroup algebra} \label{sub:radical}

We now return to the case of a LRB $S$.  Our study of the semigroup
algebra makes no use of the fact that the scalars are real numbers or
that $\{w_x\}$ is a probability distribution.  We therefore work in
the semigroup algebra $kS$ of $S$ over an arbitrary field~$k$.

The axiom \eqref{e:homom} for LRBs says that the support map $S\onto
L$ is a semigroup homomorphism, where $L$ is viewed as a semigroup
under the join operation, $X,Y\mapsto X\vee Y$.  Extending to linear
combinations, we obtain a $k$-algebra surjection
\[
\supp\colon kS \onto kL.
\]
Now Solomon~\cite{\solomonb} showed that the semigroup algebra $kL$ is
isomorphic to a product of copies of~$k$; see also \cite{\greene} and
\cite[Section 3.9]{\stanley}.  Explicitly, if $k^L$ denotes the ring
of functions from $L$ to~$k$, then there is an algebra isomorphism
\[
\phi\colon kL \isoto k^L
\]
such that $\phi(X)$ for $X\in L$ is the function $1_{Y\ge X}$, whose
value at $Y$ is 1 if $Y\ge X$ and 0 otherwise.  (Note that $\phi$
preserves products because $1_{Y\ge X}1_{Y\ge X'}=1_{Y\ge X\vee X'}$.)
Composing $\phi$ with the support map, we obtain a map $\psi\colon
kS\onto k^L$, which plays the role of the Fourier transform.  It is
not an isomorphism but, as we will see shortly, its kernel is
nilpotent; this turns out to be enough to let us compute eigenvalues.

Before proceeding to the analysis of the kernel, we record for future
reference an explicit formula for $\phi^{-1}$.  Let
$\{\delta_X\}_{X\in L}$ be the standard basis of~$k^L$; thus
$\delta_X(Y)=1_{Y=X}$.  Then $\phi$ is given by $\phi(X)= \sum_{Y\ge
X} \delta_Y$; hence $X=\sum_{Y\ge X} \phi^{-1}(\delta_Y)$, and
M\"{o}bius inversion gives $\phi^{-1}(\delta_X)=e_X$, where
\begin{equation} \label{e:priml}
e_X=\sum_{Y\ge X}\mu(X,Y) Y.
\end{equation}
The elements $e_X$ therefore give a basis of $kL$ consisting of
pairwise orthogonal idempotents, i.e., $e_X^2=e_X$ and $e_Xe_Y=0$ for
$X\ne Y$.  In the standard terminology of ring theory, they are the
\emph{primitive idempotents} of~$kL$.

Consider now the kernel $J$ of $\supp: kS\onto kL$.  It consists of
linear combinations of elements of~$S$ such that if we lump the terms
according to supports, the coefficient sum of each lump is zero.  Thus
$J=\sum_{X\in L} J_X$, where $J_X$ consists of linear combinations
$\sum_{\supp x =X} a_x x$ with $\sum_x a_x=0$.  Suppose we compute a
product $ab$ with $a=\sum a_x x\in J_X$ and $b=\sum b_y y \in J_Y$.
If $Y\le X$, we get 0, because our axiom~\eqref{e:delete}
(Section~\ref{sub:defn}) implies that $xb=0$ for each $x$ with $\supp
x=X$.  If $Y\nleq X$, on the other hand, then $ab\in J_{X\vee Y}$, and
$X\vee Y > X$.

Next, suppose we compute a product $abc\cdots$ of several factors,
coming from $J_X$, $J_Y$, $J_Z$,\dots.  By what we have just shown, we
either get 0 or we get an increasing chain $X<X\vee Y<X\vee
Y\vee Z<\cdots$.  Since $L$ is finite, we must in fact get 0 if there
are enough factors.  Thus the ideal $J$ is nilpotent.  Summarizing, we
have:

\begin{theorem}[Bidigare] \label{t:radical}
There is an algebra surjection $\psi\colon kS\onto k^L$ whose kernel
$J$ is nilpotent.  The $X$-component of $\psi$ is the homomorphism
$\chi_X\colon kS\to k$ given by
\[
\chi_X(y)=1_{\supp y\le X}
\]
for $y\in S$.
\end{theorem}

One can express the first sentence of Theorem~\ref{t:radical} by
saying that $J$ is the radical of the ring~$kS$ and that $k^L$ is the
semisimple quotient.  Standard ring theory now implies:

\begin{corollary}
Every irreducible representation of $kS$ is 1-dimensional.  There is
one such for each $X\in L$, given by the character~$\chi_X$.
\end{corollary}

We give the proof, for the convenience of readers not familiar with
the concepts of radical and semisimple quotient.

\begin{proof}
Let $V$ be an irreducible $kS$-module.  Then $JV$ is a submodule of
$V$, so it is either $V$ or 0.  (Here $JV$ is the set of finite sums
$\sum_i a_i v_i$ with $a_i\in J$ and $v_i\in V$.)  It cannot be $V$,
because then we would have $J^mV=V$ for all $m$, contradicting the
fact that $J$ is nilpotent and $V\ne0$.  So $JV=0$, and the action of
$kS$ on $V$ factors through the quotient $k^L$.  Now consider the
action on $V$ of the standard basis vectors $\delta_X$ of $k^L$.  Each
$\delta_X V$ is a submodule (because $k^L$ is commutative), so it is
either $V$ or 0.  There cannot be more than one $X$ with $\delta_XV=V$
because $\delta_X\delta_Y=0$ for $X\ne Y$.  Since $\sum_X \delta_X=1$,
it follows that exactly one $\delta_X$ is nonzero on $V$, and it acts
as the identity.  Hence every $a\in kS$ acts on $V$ as multiplication
by the scalar $\chi_X(a)$, and irreducibility now implies that $V$ is
1-dimensional.
\end{proof}

For any finite-dimensional $kS$-module $V$, we can take a
composition series 
\[
0=V_0<V_1<\cdots<V_n=V
\]
and apply the corollary to each factor $V_i/V_{i-1}$.  It follows that
there are $X_i\in L$, $i=1,\dots,n$, such that the elements $a\in kS$
are simultaneously triangularizable on~$V$, with diagonal entries
$\chi_{X_1}(a),\dots,\chi_{X_n}(a)$.  In particular, we can read off
the eigenvalues of $a$ acting on~$V$ as soon as we know, for each
$X\in L$, how many times $\chi_X$ occurs as a composition factor.

\subsection{The eigenvalues of $P$} \label{sub:eigen1}

We can now prove the formula for the eigenvalues of our transition
matrix $P$ stated in Theorem~\ref{t:main}, in somewhat greater
generality:

\begin{theorem}
Let $S$ be a finite LRB with identity, let $k$ be a field, and let
$w=\sum_{x\in S} w_x x$ be an arbitrary element of $kS$.  Let $P$ be
defined by equation~\eqref{e:transition}.  Then $P$ has an eigenvalue
\[
\lambda_X = \sum_{\supp y\le X} w_y
\]
for each $X\in L$, with multiplicity $m_X$, where
\begin{equation} \label{e:mult1}
\sum_{Y\ge X} m_Y = c_X
\end{equation}
for each $X\in L$.
\end{theorem}

\begin{proof}
Recall from Section~\ref{sub:algebra} that the eigenvalues
of~$P$ are the same as the eigenvalues of $w$ acting by left
multiplication on the ideal $kC\subseteq kS$.  For each $X\in L$, let
$m'_X$ be the number of composition factors of $kC$ given by the
character $\chi_X$.  Then the discussion at the end of
Section~\ref{sub:radical} shows that $P$ has eigenvalues $\chi_X(w)$
with multiplicity~$m'_X$.  Now
\[
\chi_X(w)=\sum_{y\in S} w_y 1_{\supp y\le X} = \sum_{\supp y \le X} w_y=\lambda_X,
\]
so the proof will be complete if we show that $\sum_{Y\ge X} m'_Y =
c_X$ for all $X\in L$.

Consider an arbitrary $x\in S$.  It acts on $kC$ as an idempotent
operator, projecting $kC$ onto the linear span of the chambers in
$S_{\ge x}$.  The rank~$r$ of this projection is therefore the number
$c_X$ defined in Section~\ref{s:main}, where $X=\supp x$.  On the
other hand, the rank of a projection is the multiplicity of 1 as an
eigenvalue, so
\[
r=\sum_{\substack{Y\in L\\ \chi_Y(x)=1}} m'_Y = \sum_{Y\ge X} m'_Y.
\]
Equating the two expressions for $r$ gives $\sum_{Y\ge X}
m'_Y=c_X$, as required.
\end{proof}

We turn now to the proof that $P$ is diagonalizable when $k=\R$ and
$w$ is a probability distribution.

\section{Semisimplicity} \label{s:semisimple}

Let $\R[w]\subseteq \R S$ be the subalgebra (with identity) generated
by $w=\sum_{x\in S} w_x x$, where $w_x\ge0$ and $\sum_x w_x=1$.  We
will show that $\R[w]$ is semisimple; more precisely, it is isomorphic
to a direct product of copies of~$\R$.  This implies that the action
of~$w$ is diagonalizable in every $\R S$-module; in particular, it
implies that the transition matrix~$P$ of our walk on chambers is
diagonalizable, as asserted in Theorem~\ref{t:main}.

In order to show the idea of the proof in its simplest form, we begin
by giving in Section~\ref{sub:diag} a criterion (probably known) for
the diagonalizability of a matrix~$A$, involving the poles of the
generating function for the powers of~$A$.  In
Section~\ref{sub:criterion} we essentially repeat the proof, but in a
more abstract setting; the result is a criterion for semisimplicity of
an algebra generated by a single element~$a$, involving the generating
function for the powers of~$a$.  Then in Section~\ref{sub:powers} we
compute the powers of our element $w\in\R[w]\subseteq \R S$, and we
deduce a formula for the generating function.  The criterion of
Section~\ref{sub:criterion} is visibly satisfied, and we get the
desired semisimplicity result in Section~\ref{sub:proof}.  As a
byproduct of the proof we obtain formulas for the primitive
idempotents of~$\R[w]$, which we state in Section~\ref{sub:primitive}.
As a simple example, we write out the formulas for the Tsetlin library
with uniform weights in Section~\ref{sub:tsetlinu}.  In a very
technical Section~\ref{sub:primitive2} we attempt to organize the
formulas in a sensible way.  Finally, we return to the Tsetlin library
in Section~\ref{sub:example}, this time with arbitrary weights, to
illustrate the results of Section~\ref{sub:primitive2}.

\subsection{Diagonalizability} \label{sub:diag}

Let $M_n(\CC)$ be the ring of $n\times n$ matrices over $\CC$, and let
$A\in M_n(\CC)$.  [With minor changes we could work over an arbitrary
field instead of~$\CC$.]  Consider the generating function
\[
f(t)=\sum_{m\ge0} A^m t^m = \frac{1}{I-tA}\,,
\]
where $I$ is the identity matrix and the fraction is to be interpreted
as $(I-tA)^{-1}$.  The series converges for small $t\in\CC$ and
represents a holomorphic function with values in $M_n(\CC)$.  It is
initially defined in a neighborhood of~0, but we will see in
Proposition~\ref{p:diag} that $f$ is a rational function, i.e., that
each of the $n^2$ matrix entries is a rational function in the usual
sense.  Let
\[
g(z)= (1/z)f(1/z)=\frac{1}{zI-A}\,,
\]
initially defined for $z$ in a neighborhood of $\infty$.

\begin{proposition} \label{p:diag}
The function $g$ is rational, with poles precisely at the eigenvalues
of~$A$.  The matrix $A$ is diagonalizable if and only if the poles
of~$g$ are all simple.  In this case $g$ has a partial fractions
decomposition of the form
\[
g(z)=\sum_i\frac{E_i}{z-\lambda_i},
\]
where the $\lambda_i$ are the distinct eigenvalues of~$A$ and $E_i$ is
the projection onto the $\lambda_i$-eigenspace.
\end{proposition}

``Projection'' here refers to the decomposition of $\CC^n$ into
eigenspaces.

\begin{proof}
Consider the Jordan decomposition $A=\sum_i (\lambda_iE_i + B_i)$;
here the $E_i$ are pairwise orthogonal idempotents summing to~$I$, the
$B_i$ are nilpotent, and $B_i=B_iE_i=E_iB_i$.  If $A$ is
diagonalizable, then each $B_i=0$ and we have
\[
g(z)=\frac{1}{zI-A}=\sum_i\frac{1} {z-\lambda_i}E_i,
\]
as required.  If $A$ is not diagonalizable, then for some eigenvalue
$\lambda_i$ we have $B_i\ne0$.  Since $g(z)$ can be computed in each
Jordan block separately, we may assume that $A=\lambda I + B$, where
$B^r=0$ but $B^{r-1}\ne0$ for some $r>1$.  Then
\begin{align*}
g(z)&=\frac{1}{zI-A}\\
&=\frac{1}{(z-\lambda)I-B}\\
&=\frac{1}{z-\lambda}\cdot\frac{1}{I-(z-\lambda)^{-1}B}\\
&=\sum_{j=0}^{r-1} \frac{B^j}{(z-\lambda)^{j+1}}\,.
\end{align*}
Thus $g(z)$ is rational and has a pole of order $r>1$ at $z=\lambda$.
\end{proof}

\subsection{A semisimplicity criterion} \label{sub:criterion}

The ring-theoretic version of what we have just done goes as follows.
Let $k$ be a field and $R$ a finite-dimensional commutative
$k$-algebra (with identity).  For simplicity, we will pretend that $k$
is a subfield of~$\CC$, so that we can speak of convergent power
series; to deal with a general field~$k$, one needs to work with formal
power series.  In our application we will have $k=\R$.

Assume that $R$ is generated by a single element~$a$.  Thus $R\iso
k[x]/(p)$ for some polynomial $p$, with $a$ corresponding to
$x$~mod~$p$.  We give here a criterion for $R$ to be \emph{split
semisimple}, i.e., isomorphic to $k^I$, a product of copies of $k$
indexed by a (finite) set $I$.  Giving such an isomorphism is
equivalent to giving a basis $(e_i)_{i\in I}$ for $R$ consisting of
pairwise orthogonal idempotents.  The $e_i$ are then characterized as
the \emph{primitive} idempotents of $R$, i.e., the nonzero idempotents
that cannot be decomposed as sums of pairwise orthogonal nonzero
idempotents.

Consider the generating function
\[
f(t)=\sum_{m=0}^\infty a^m t^m =\frac{1}{1-at}\,,
\]
where the fraction is to be interpreted as $(1_R-at)^{-1}$.  It will
follow from the proof of Proposition~\ref{p:criterion} that the series
has a positive radius of convergence and that $f$ is a rational
function with values in~$A$; this means that if we express $f(t)$ in
terms of a basis for~$A$, then each component is a rational function
in the usual sense.  Let
\[
g(z)= (1/z)f(1/z)=\frac{1}{z-a}\,;
\]
here we identify $k$ with the ring of scalar multiples of the
identity~$1_R$, so that $z-a$ means $z1_R -a$.

\begin{proposition} \label{p:criterion}
The $k$-algebra $R$ is split semisimple if and only if $g(z)$ has the
form
\begin{equation} \label{e:gen}
g(z)= \sum_{i\in I} \frac{e_i}{z-\lambda_i},
\end{equation}
where the $\lambda_i$ are distinct elements of $k$ and the $e_i$ are
nonzero elements of $R$.  In this case the $e_i$ are the primitive
idempotents of $R$, and the generator $a$ of $R$ is equal to
$\sum_{i\in I} \lambda_i e_i$.
\end{proposition}

\begin{proof}
Suppose $A$ is split semisimple with primitive idempotents
$(e_i)_{i\in I}$, and write $a=\sum_i \lambda_i e_i$.  Then
$a^m=\sum_i \lambda_i^m e_i$, $f(t)=\sum_i (1-\lambda_it)^{-1} e_i$,
and the expression \eqref{e:gen} for $g(z)=(1/z)f(1/z)$ follows at once.

Conversely, suppose $A$ is not split semisimple.  Assume first that
the minimal polynomial $p$ of $a$ splits into linear factors in
$k[x]$, say $p(x)=\prod_{i\in I} (x-\lambda_i)^{r_i}$, where the
$\lambda_i$ are distinct.  By the Chinese remainder theorem,
\begin{equation} \label{e:crt}
A\iso\prod_{i\in I} k[x]/(x-\lambda_i)^{r_i},
\end{equation}
and the assumption that $A$ is not split semisimple implies that some
$r_i>1$.

Since $g(z)$ can be computed componentwise with respect to the
decomposition~\eqref{e:crt} of~$A$, we may assume that there is only
one factor, i.e., that $A=k[x]/(x-\lambda)^r$ for some~$\lambda$,
where $r>1$.  Then $a=\lambda+b$, where $b^r=0$ but $b^{r-1}\ne0$;
hence
\begin{align*}
g(z)&=\frac{1}{z-a}\\
&=\frac{1}{(z-\lambda)-b}\\
&=\frac{1}{z-\lambda}\cdot\frac{1}{1-(z-\lambda)^{-1}b}\\
&=\sum_{j=0}^{r-1} \frac{b^j}{(z-\lambda)^{j+1}}\,.
\end{align*}
Thus $g(z)$ has a pole of order $r>1$ at $z=\lambda$ and hence does
not have the form~\eqref{e:gen}.

If $p$ does not split into linear factors, extend scalars to a
splitting field $k'$ of $p$ and apply the results above to
$A'=k'\otimes_k A \iso k'[x]/(p)$.  Then $g(z)$, viewed now as a
function $k'\to A'$, has poles at the roots of~$p$, at least one of
which is not in~$k$.  Once again, $g(z)$ does not have the
form~\eqref{e:gen}.
\end{proof}

\subsection{A formula for $w^m$} \label{sub:powers}

Let $S$ be a LRB and let $w=\sum_{x\in S} w_x x$, where $\{w_x\}$ is a
probability distribution on~$S$.  From now on we identify $w$ with
$\{w_x\}$ and simply say that $w$ is a probability distribution.  We
wish to apply Proposition \ref{p:criterion} to $R=\R[w]\subseteq\R S$.
To this end we need a formula for~$w^m$.  As an aid to the intuition,
we  use probabilistic language in deriving this formula.  The
interested reader can recast the discussion in purely algebraic
language, where it is valid with $\R$ replaced by an arbitrary field
$k$ and $w$ by an arbitrary element of the semigroup algebra $kS$.
Our methods in this section are inspired by the paper of
Fill~\cite{\Fill}.

By a \emph{reduced word} we mean an $l$-tuple
$\mathbf{x}=(x_1,\dots,x_l)$, $x_i\in S$, such that for each
$i=1,\dots,l$ we have $\supp x_i \nleq \supp(x_1\cdots x_{i-1})$.
Equivalently, if we set $X_i=\supp(x_1\cdots x_i)$, then we get a
strictly increasing chain 
\[
\zerohat=X_0<X_1<\cdots<X_l
\]
in~$L$.  We say that $\mathbf{x}$ is a \emph{reduced decomposition} of
the element $\bar{\mathbf{x}}=x_1x_2\cdots x_l\in S$.  The intuitive
meaning of this is that there is no obvious way to shorten the
expression $x_1x_2\cdots x_l$ by using the axiom~\eqref{e:delete} to
delete factors.  If an $m$-tuple $(x_1,\dots,x_m)$ is not necessarily
reduced, there is a reduced word $(x_1,\dots,x_m)\sphat\,$, obtained
by deleting any $x_i$ such that $\supp x_i\le \supp(x_1\cdots
x_{i-1})$.

\begin{remark}
It might seem more natural to require the ``letters'' $x_i$ in a
reduced word to be in some given generating set $S_1\subseteq S$.  In
practice, one is typically interested in $S_1=\{x\in S:w_x\ne0\}$.
Our convention of allowing arbitrary $x_i$ is harmless, however, since
only those words whose letters are in $S_1$ make a nonzero
contribution to the formula~\eqref{e:formula} that we are going to
derive.
\end{remark}

We need some notation in order to state the formula.  Let
$\mathbf{x}=(x_1,\dots,x_l)$ be a reduced word of length
$l=l(\mathbf{x})$, with associated chain
\[
\zerohat=X_0<X_1<\cdots<X_l.
\]
Let $\lambda_0,\lambda_1,\dots,\lambda_l$ be the corresponding eigenvalues
$\lambda_{X_i}$ as in Theorem~\ref{t:main}, and, for $n\ge0$, let
$H_n(\mathbf{x})=h_n(\lambda_0,\dots,\lambda_l)$, where $h_n$ is the
complete homogeneous symmetric function of degree~$n$ (sum of all
monomials of degree~$n$).  Let $w_{\mathbf{x}}=w_{x_1}w_{x_2}\cdots
w_{x_l}$.

\begin{proposition}
Let $S$ be a finite LRB with identity and let $w\in\R S$ be a
probability distribution.  For any $m\ge0$,
\begin{equation} \label{e:formula}
w^m =\sum_{\mathbf{x}} H_{m-l(\mathbf{x})}(\mathbf{x})w_{\mathbf{x}} \bar{\mathbf{x}},
\end{equation}
where $\mathbf{x}$ ranges over the reduced words of length
$l(\mathbf{x})\le m$.
\end{proposition}

\begin{proof}
Let $y_1,y_2,\cdots, y_m$ be independent picks from the probability
measure~$w$.  We will get a formula for $w^m$ by computing the
distribution of the reduced word $(y_1,\dots,y_m)\sphat\,$; for we have
\begin{equation} \label{e:distribution}
w^m = \sum_{l(\mathbf{x})\le m} \Pr\{(y_1,\dots,y_m)\sphat=\mathbf{x}\}\bar{\mathbf{x}}.
\end{equation}
Given a reduced word $\mathbf{x}=(x_1,\dots,x_l)$ with associated
chain $(X_0,\dots,X_l)$, we compute the probability
in~\eqref{e:distribution} as follows.  Let $S_i=S_{\le X_i} = \{x\in S:
\supp x\le X_i\}$, and let $\lambda_i=\lambda_{X_i}$.  In order to
have $(y_1,\dots,y_m)\sphat = \mathbf{x}$, the $m$-tuple
$(y_1,\dots,y_m)$ must consist of $i_0$ elements of~$S_0$, then $x_1$,
then $i_1$ elements of~$S_1$, then $x_2$, and so on, ending with $i_l$
elements of~$S_l$, where $i_0,\dots,i_l\ge0$ and $i_0+\cdots+i_l=m-l$.
The probability of this, for fixed $i_0,\dots,i_l$, is
$\lambda_0^{i_0}w_{x_1}\lambda_1^{i_1}w_{x_2}\lambda_2^{i_2}\cdots
w_{x_l} \lambda_l^{i_l}$.  Summing over all possible
$(i_0,\dots,i_l)$, we see that the probability in question is
$H_{m-l}(\mathbf{x})w_{\mathbf{x}}$, whence~\eqref{e:formula}.
\end{proof}

Formula \eqref{e:formula}  can be rewritten in terms of the
function $g(z)$ of Section~\ref{sub:criterion}.  Given a reduced word
$\mathbf{x}$ as above, set
\[
g_{\mathbf{x}}(z) = \prod_{i=0}^l \frac{1}{(z-\lambda_i)}.
\]

\begin{corollary}
Let $g(z)=(1/z)f(1/z)$, where $f(t)=\sum_{m\ge0} w^m t^m$.  Then
\begin{equation} \label{e:gen1}
g(z)=\sum_{\mathbf{x}} g_{\mathbf{x}}(z) w_{\mathbf{x}} \bar{\mathbf{x}},
\end{equation}
where $\mathbf{x}$ ranges over all reduced words.
\end{corollary}

\begin{proof}
Fix a reduced word $\mathbf{x}$ of length $l$, and let
$\lambda_0,\dots,\lambda_l$ be as above.  Then
\begin{align*}
\sum_{m\ge l} H_{m-l}(\mathbf{x})t^m &=
t^l\sum_{m\ge0}h_m(\lambda_0,\dots,\lambda_l)t^m\\
&= t^l \prod_{i=0}^l \frac{1}{1-\lambda_i t}.
\end{align*}
Setting $t=1/z$ and multiplying by $1/z$, we obtain $g_{\mathbf{x}}(z)$;
\eqref{e:gen1} now follows from~\eqref{e:formula}.
\end{proof}

\subsection{Proof of semisimplicity} \label{sub:proof}

Call an element $X\in L$ \emph{feasible} for $w$ if $X=\supp(x_1\cdots
x_m)$ with $w_{x_i}\ne0$ for $i=1,\dots,m$ or, equivalently, if $X$ is
the join of elements $\supp x$ with $w_x\ne0$.  Let $L_w$ be the set
of feasible elements of~$L$.  In applying formulas \eqref{e:formula}
and \eqref{e:gen1}, we need only consider reduced words~$\mathbf{x}$
whose associated chain is in~$L_w$, since otherwise
$w_{\mathbf{x}}=0$.  The eigenvalues $\lambda_0,\dots,\lambda_l$ are
then all distinct; in fact, we have
$\lambda_0<\lambda_1<\cdots<\lambda_l$.  So we obtain an expression of
the form~\eqref{e:gen} for $g(z)$ by splitting each
$g_{\mathbf{x}}(z)$ into partial fractions.  We have therefore proved
the first part of the following theorem:

\begin{theorem}
Let $S$ be a finite LRB with identity and let $w\in\R S$ be a
probability distribution.  Then the subalgebra $\R[w]$ is split
semisimple.  Consequently, the action of~$w$ on any $\R S$-module is
diagonalizable.
\end{theorem}

The second assertion is an easy consequence of the first.  Indeed, if
we write $w=\sum_i \lambda_i e_i$, where the $e_i$ are the primitive
idempotents of~$\R[w]$, then any $\R S$-module $V$ decomposes as
$V=\bigoplus_i e_iV$, with $w$ acting as multiplication by~$\lambda_i$
on~$e_iV$.

Theorem 1 is now completely proved.

\begin{remark}
Everything we have done remains valid with $\R$ replaced by an
arbitrary field~$k$ and $w$ by an arbitrary element of~$kS$, with one
proviso.  Namely, it is no longer automatic that the eigenvalues
$\lambda_0,\dots,\lambda_l$ are distinct.  In order to guarantee this,
we need to assume that $w$ satisfies the following condition: Whenever
$X<Y$ in~$L_w$, one has $\lambda_X \ne \lambda_Y$.  Under this
assumption, then, $k[w]$ is split semisimple.
\end{remark}

\subsection{Primitive idempotents, first version} \label{sub:primitive}

It is easy to determine the primitive idempotents of $\R[w]$ (or
$k[w]$, under the hypotheses of the remark above) by using
Proposition~\ref{p:criterion} and formula~\eqref{e:gen1}.  We assume,
without loss of generality, that $S$ is generated by $\{x\in
S:w_x\ne0\}$; this implies that $L_w=L$.

Suppose first that $w$ is \emph{generic}, by which we mean that
$\lambda_X\ne\lambda_Y$ for $X\ne Y$ in $L$.  (Thus we are excluding
those probability measures that lie on the union of a certain finite
collection of hyperplanes in~$\R S$.)  Then the homomorphism
$\psi\colon\R S\onto\R^L$ of Theorem~\ref{t:radical}
(Section~\ref{sub:radical}) maps $\R[w]$ onto $\R^L$; in fact,
$\psi(w)=\sum_{X\in L} \lambda_X \delta_X$, and it is easy to check
that $\R^L$ is generated as an algebra by any element whose components
are all distinct.  Since $\R[w]$ is known to be semisimple and $\ker
\psi$ is nilpotent, it follows that $\psi$ maps $\R[w]$ isomorphically
onto $\R^L$.  Hence $\R[w]$ has one primitive idempotent $e_X$ for
each $X\in L$, and
\begin{equation} \label{e:decomp}
w=\sum_{X\in L} \lambda_X e_X.
\end{equation}
To compute $e_X$, we have to multiply the right side of \eqref{e:gen1}
by $z-\lambda_X$ and then set $z=\lambda_X$.

Let $\mathbf{x}$ be a reduced word as in Section~\ref{sub:powers}, and
suppose its associated chain passes through $X$, say $X=X_i$.  Then
the residue of $g_{\mathbf{x}}(z)$ at $z=\lambda_X$ is
\[
R_{X,\mathbf{x}} = \frac{(-1)^{l-i}}{(\lambda_X-\lambda_0)\cdots
(\lambda_X-\lambda_{i-1}) (\lambda_{i+1}-\lambda_X)\cdots
(\lambda_l-\lambda_X)}\,.
\]
Hence
\begin{equation} \label{e:formulai}
e_X=\sum_{\mathbf{x}} R_{X,\mathbf{x}} w_{\mathbf{x}} \bar{\mathbf{x}},
\end{equation}
where $\mathbf{x}$ ranges over all reduced words whose chain passes
through~$X$.

If $w$ is not generic, then this formula still makes sense, but one
has to sum the $e_X$ having a common value of $\lambda_X$ in order to
get the primitive idempotents of~$\R[w]$; the $e_X$ themselves may not
lie in~$\R[w]$.  Notice, however, that the $e_X$ still form an
orthogonal family of idempotents in~$\R S$ summing to~1, and the
decomposition of~$w$ given in~\eqref{e:decomp} is still valid.  To see
this, note that these assertions can be formulated as polynomial
equations in the variables~$w_x$; since the equations are valid
generically, they must hold as algebraic identities.

Summarizing, we have:

\begin{theorem}
Let $S$ be a finite LRB with identity and let $w$ be a probability
distribution on $S$.  Assume that $S$ is generated, as a semigroup
with identity, by $\{x\in S:w_x\ne0\}$.  Then~\eqref{e:formulai}
defines an orthogonal family of idempotents~$e_X$, $X\in L$, such that
\eqref{e:decomp} holds.  If $w$ is generic, then the $e_X$ are the
primitive idempotents of~$\R[w]$.  In general, the decomposition of
$w$ as a linear combination of primitive idempotents of~$\R[w]$ is
obtained by grouping the terms in~\eqref{e:decomp} according to the
value of~$\lambda_X$.
\end{theorem}

\subsection{Example: The Tsetlin library with uniform weights}
\label{sub:tsetlinu}

Let $S=F_n$, with uniform weights $w_i=1/n$ on the
elements $(i)$ of length~1.  For each $x=(x_1,\dots,x_l)\in S$, the
only reduced decomposition $\mathbf{x}$ of~$x$ with
$w_{\mathbf{x}}\ne0$ is the obvious one,
$\mathbf{x}=((x_1),\dots,(x_l))$.  The associated chain is given by
$X_i=\{x_1,\dots,x_i\}$ for $0\le i\le l$.  If $X=X_i$, then the
contribution of~$\mathbf{x}$ to~$e_X$ is
\begin{equation} \label{e:formulaiu}
R_{X,\mathbf{x}}w_{\mathbf{x}}x = (-1)^{l-i}\frac{x}{i!(l-i)!}.
\end{equation}
To get $e_X$, then, we have to sum over all $x$ having some ordering
of~$X$ as an initial segment.  The eigenvalue corresponding to~$X$ is
$\lambda_X=i/n$.  We conclude that $\R[w]$ has $n+1$ primitive
idempotents $e_0,e_1,\dots,e_n$, where $e_i=\sum_{\abs{X}=i} e_X$;
hence $e_i$ is obtained by summing the right-hand side
of~\eqref{e:formulaiu} over all $x$ of length $l\ge i$.  If
$\sigma_l\in\R S$ is the sum of all $x\in S$ of length $l$, the result is
\begin{equation*} 
e_i = \sum_{l=i}^n
(-1)^{l-i}\frac{\sigma_l}{i!(l-i)!}= \sum_{l=i}^n (-1)^{l-i}
\binom{l}{i}\frac{\sigma_l}{l!}\,.
\end{equation*}
The decomposition of $w$ is
\[
w=\frac{1}{n}\sum_{i=0}^n i e_i.
\]

Recall that the Tsetlin library can also be obtained by using a
quotient $\bar{S}=\bar{F}_n$ of~$F_n$ (Section~\ref{sub:braidsemi}).
For any $a\in\R S$, let $\bar{a}$ be its image in~$\R\bar{S}$.  Then
the probability distribution on~$\bar{S}$ that gives the Tsetlin
library with uniform weights is~$\bar{w}$.  One can check that the
quotient map $S\onto\bar{S}$ induces a surjection
$\R[w]\onto\R[\bar{w}]$ with 1-dimensional kernel, spanned
by~$e_{n-1}$.  Thus $\R[\bar{w}]$ has $n$ primitive idempotents
$\bar{e}_0,\dots,\bar{e}_{n-2},\bar{e}_n$, with
\begin{equation} \label{e:tsetlinu}
\bar{e}_i = \sum_{l=i}^n (-1)^{l-i}
\binom{l}{i}\frac{\bar{\sigma}_l}{l!}\,.
\end{equation}
This equation is also valid for $i=n-1$, in which case its content is
that $\bar{e}_{n-1}=0$, as stated above; this follows from the fact
that $\bar{\sigma}_{n-1} = \bar{\sigma}_n$.  Formula
\eqref{e:tsetlinu} is essentially the same as a formula of
Diaconis--Fill--Pitman~\cite[(4.5)]{\dfp}, except that these authors
work with operators on $\R C$ and interpret the answer in terms of
Solomon's descent algebra.  We will explain this in more detail in
Section~\ref{sub:descentwalk}.

\begin{remark}
We could equally well have treated general weights, but instead we
will do that in Section~\ref{sub:example}, as an illustration of
a different version of the formula for~$e_X$.
\end{remark}

\subsection{Primitive idempotents, second version} \label{sub:primitive2}

In this quite technical subsection we attempt to make sense out of
formula~\eqref{e:formulai} for the primitive idempotents.  Our goal,
motivated by equation~\eqref{e:priml} for the primitive idempotents
in~$\R L$, is to write~\eqref{e:formulai} in the form
\begin{equation} \label{e:formula1}
e_X=\sum_{Y\ge X}\nu_{X,Y}\,,
\end{equation}
where $\nu_{X,Y}$ is a certain signed measure on $C_Y=\{y\in S: \supp
y=Y\}$.  Comparing this with~\eqref{e:priml}, we see that $\nu_{X,Y}$
necessarily has total mass $\mu(X,Y)$.  As usual, $\nu_{X,Y}$ is
identified with a linear combination of the elements of $C_Y$, hence
it is an element of $\R S$ and \eqref{e:formula1} makes sense.  The
definition of $\nu_{X,Y}$ is complicated.  We begin with the case
$X=\zerohat$, which is slightly simpler.

Fix  $Y\in L$ and consider an arbitrary chain $\mathbf{X}$ from
$\zerohat$ to $Y$,
\[
\zerohat=X_0<X_1<\cdots <X_l=Y,
\]
of length $l=l(\mathbf{X})\ge0$.  We associate to $\mathbf{X}$ a
defective probability measure $\rho_{\mathbf{X}}$ on $C_Y$, as
follows.  Given $y\in C_Y$, consider all reduced decompositions
$(y_1,\dots,y_l)$ of $y$ whose associated chain is $\mathbf{X}$.  For
each such decomposition, set $w_i=w_{y_i}$ and
$\lambda_i=\lambda_{X_i}$ and form the product
\[
\frac{w_1}{\lambda_1-\lambda_0} \frac{w_2}{\lambda_2-\lambda_0} \cdots
\frac{w_l}{\lambda_l-\lambda_0}.
\]
Then $\rho_{\mathbf{X}}(y)$ is the sum of all these products.  This
has a probabilistic interpretation: Pick elements $y_i\in S_{\le
X_i}-\{\text{id}\}$ independently, according to the weights $w_y$,
where $i=1,\dots,l$.  If $\supp(y_1\cdots y_i)=X_i$ for each $i$, form
the product $y_1\cdots y_l$.  This defines a defective random variable
with values in $C_Y$, and $\rho_{\mathbf{X}}$ is its distribution.  (A
defective random variable is one that is defined with probability
$\le1$; its distribution is a positive measure having total mass
$\le1$.)  The signed measure $\nu_{\zerohat,Y}$ is now obtained by
taking an alternating sum:
\[
\nu_{\zerohat,Y} = \sum_{\mathbf{X}} (-1)^{l(\mathbf{X})}\rho_{\mathbf{X}},
\]
where $\mathbf{X}$ ranges over all chains from $\zerohat$ to $Y$.

For general $X$, consider chains $\mathbf{X}$ from $X$ to $Y$:
\[
X=X_0<X_1<\cdots <X_l=Y.
\]
We first define a defective probability measure $\rho_{x,\mathbf{X}}$
on $C_Y$, depending on a choice of $x\in C_X$.  Pick $y_1,\dots,y_l$
independently, with $y_i\in S_{\le X_i}-S_{\le X}$.  If
$\supp(xy_1\cdots y_i)=X_i$ for each $i$, form the product $xy_1\cdots
y_l$.  This gives a defective random variable with values in $C_Y$,
and $\rho_{x,\mathbf{X}}$ is its distribution.  We now set
\begin{equation} \label{e:nux}
\nu_{x,Y} = \sum_{\mathbf{X}} (-1)^{l(\mathbf{X})}\rho_{x,\mathbf{X}}\,,
\end{equation}
where $\mathbf{X}$ ranges over all chains from $X$ to $Y$.  We can
also describe $\nu_{x,Y}$ as the measure obtained by applying the
procedure of the previous paragraph to $S_{\ge x}$, using the
probability measure obtained from $w$ via the projection
$S\onto S_{\ge x}$.

We now define the desired $\nu_{X,Y}$ by averaging over~$x\in C_X$:
\begin{equation} \label{e:nuX}
\nu_{X,Y} = \sum_{x\in C_X} \pi_X(x)\nu_{x,Y}\,,
\end{equation}
where $\pi_X$ is the stationary distribution of the random walk on
$C_X$ driven by the weights $w_y$, $y\in S_{\le X}$ (scaled to give a
probability distribution).  This completes the formula for~$e_X$.  We
leave it to the interested reader to verify that \eqref{e:formula1} is
indeed a reformulation of~\eqref{e:formulai}; the starting point is to
group the terms in~\eqref{e:formulai} according to the chain
$\mathbf{X}$ associated with~$\mathbf{x}$.

\begin{remark}
The idempotent $e_{\onehat}$ is the stationary distribution $\pi$ of
our random walk.  The decomposition of $w$ can therefore be written as
\[
w=\pi+\sum_{X<\onehat} \lambda_X e_X,
\]
so that
\[
w^m = \pi + \sum_{X<\onehat}\lambda_X^m e_X.
\]
In theory, this should make it possible to give precise estimates for
\[
\norm{P_c^m -\pi} = \frac{1}{2}\Bigl\|\sum_{X<\onehat}\lambda_X^m e_Xc\Bigr\|_1.
\]
In practice, however, the presence of signs makes this very tricky.
\end{remark}

\subsection{Example: The Tsetlin library} \label{sub:example}

We return to $S=F_n$ and the Tsetlin library, but now with generic
weights $w_1,\dots,w_n$.  Since the associated lattice $L$ is the
Boolean lattice of subsets of~$[n]$, we get $2^n$ primitive
idempotents $e_X$ in~$\R[w]$, and $w=\sum_{X\subseteq[n]} \lambda_X
e_X$, with $\lambda_X=\sum_{x\in X} w_x$.  Working through the
definition of the signed measure $\nu_{X,Y}$ for $X\subseteq Y$, one
finds that it is $(-1)^{\abs{Y-X}}$ times the distribution of the
following random ordering of $Y$: Sample without replacement from $X$,
getting an ordering $(x_1,\dots,x_i)$ of~$X$, where $i=\abs{X}$;
sample without replacement from $Y-X$, getting an ordering
$(y_1,\dots,y_j)$ of $Y-X$, where $j=\abs{Y-X}$; now form
$(x_1,x_2,\dots,x_i,y_j,\dots,y_2,y_1)$.  Note the reversal of the
ordering of the $y$'s; thus we are building a random ordering of~$Y$
by accumulating elements of~$X$ from left to right and elements of
$Y-X$ from right to left.

This gives a very explicit formula for the primitive idempotent
\[
e_X=\sum_{Y\supseteq X}\nu_{X,Y}.
\]
This and the related formula
\[
w^m = \sum_{X\subseteq[n]} \lambda_X^m e_X,
\]
are essentially the formulas of Fill~\cite{\Fill}, except that he
works with operators on $\R C$ instead of with elements of $\R S$.  To
get his formulas, right multiply the formulas above by a permutation
$\sigma$ (chamber of $S$) and pick out the $\tau$-component.  This
gives $(\sigma,\tau)$-entries of matrices.  One can check that $e_X$
annihilates $\R C$ if $\abs{X}=n-1$, so the spectral decomposition of
left multiplication by~$w$ on~$\R C$ only involves $2^n-n$
idempotents, as in Fill's paper.  This would have arisen more
naturally if we had used $\bar{F}_n$ instead of~$F_n$.

The reader who has come this far might find it a useful exercise to
rederive the formulas for the uniform case
(Section~\ref{sub:tsetlinu}) from those above.

\section{Reflection arrangements and Solomon's descent algebra} \label{s:descent}

The work of Bayer--Diaconis \cite{\riffle} and Diaconis--Fill--Pitman
\cite{\dfp} relates certain card-shuffling random walks on the
symmetric group $S_n$ to subalgebras of Solomon's descent
algebra~\cite{\solomon}.  We show here how this surprising connection
arises naturally from semigroup considerations.  We work with an
arbitrary finite Coxeter group~$W$ and its associated hyperplane face
semigroup $\Sigma$ (the Coxeter complex of~$W$).  But we will try to
explain everything in concrete terms for the case $W=S_n$, in an
effort to make the discussion accessible to readers unfamiliar with
Coxeter groups.

Our treatment can be viewed as an elaboration of Tits's appendix to
Solomon's paper, with further ideas borrowed from Bidigare's
thesis~\cite{\bid}.  In particular, we  use (and include a proof
of) Bidigare's theorem that Solomon's descent algebra is
anti-isomorphic to the $W$-invariant part of the semigroup algebra
of~$\Sigma$.

In this section the probability measure driving our random walk on the
chambers of $\Sigma$ is denoted by~$p$ instead of $w$, so that we
can reserve the letter $w$ for a typical element of $W$.

\subsection{Finite reflection groups} \label{sub:reflection}

We begin with a very quick review of the basic facts that we need
about finite Coxeter groups and their associated simplicial complexes
$\Sigma$.  Details can be found in many places, such as
\cite{\brown,\gb,\humphreys,\tits}.  A \emph{finite reflection group}
on a real inner-product space $V$ is a finite group of orthogonal
transformations of $V$ generated by reflections $s_H$ with respect to
hyperplanes $H$ through the origin.  The set of hyperplanes $H$ such
that $s_H\in W$ is the \emph{reflection arrangement} associated
with~$W$.  Its hyperplane face semigroup $\Sigma$ can be identified
with the set of simplices of a simplicial complex, called the
\emph{Coxeter complex} of~$W$.  Geometrically, this complex is gotten
by cutting the unit sphere in~$V$ by the hyperplanes~$H$, as in
Section~\ref{sub:cell}.  (As explained there, one might have to first
pass to a quotient of~$V$.)  The action of~$W$ on~$V$ induces an action of
$W$ on~$\Sigma$, and this action is simply-transitive on the chambers.
Thus the set $\C$ of chambers can be identified with~$W$, once a
``fundamental chamber'' $C$ is chosen.

The canonical example is $W=S_n$, acting on $\R^n$ by permuting the
coordinates.  The arrangement in this case is the braid arrangement
(Section~\ref{sub:braid}).  The Coxeter complex $\Sigma$ can be
identified with the following abstract simplicial complex:  The
vertices are the proper nonempty subsets $X\subset[n]=\{1,\dots,n\}$,
and the simplices are the chains of such subsets.  The $S_n$-action is
induced by the action of $S_n$ on~$[n]$.  The product on~$\Sigma$ was
discussed in Section~\ref{sub:braidsemi}, where $\Sigma$ was
identified with the semigroup~$\B$ of ordered partitions.   The
chambers of~$\Sigma$ correspond to permutations $w$ of $[n]$, with $w$
corresponding to the chamber
\[
\{w(1)\}<\{w(1),w(2)\}<\cdots<\{w(1),w(2),\dots,w(n-1)\}.
\]
This is the same as the identification of $\C$ with $S_n$ that results
from choosing
\[
\{1\}<\{1,2\}<\cdots<\{1,2,\dots,n-1\}
\]
as fundamental chamber.

\subsection{Types of simplices}

The number $r$ of vertices of a chamber of $\Sigma$ is called the
\emph{rank} of $\Sigma$ (and of $W$); thus the dimension of $\Sigma$
as a simplicial complex is $r-1$.  It is known that one can color the
vertices of $\Sigma$ with $r$ colors in such a way that vertices
connected by an edge have distinct colors.  The color of a vertex is
also called its \emph{label}, or its \emph{type}, and we denote by $I$
the set of all types.  We can also define $\type(F)$ for any
$F\in\Sigma$; it is the subset of~$I$ consisting of the types of the
vertices of~$F$.  For example, every chamber has type $I$, while the
empty simplex has type $\emptyset$.  The action of $W$ is
type-preserving; moreover, two simplices are in the same $W$-orbit if
and only if they have the same type.  In our canonical example with
$W=S_n$, the rank is $n-1$, the set of types is $I=\{1,\dots,n-1\}$,
and $\type(X)=\abs{X}$ for any vertex~$X$ (proper nonempty subset
of~$[n]$).

The labeling allows us to refine the adjacency relation on chambers
defined in Section~\ref{sub:gallery}:  If $C,C'$ are distinct adjacent
chambers and $F$ is their common face of codimension~1, then
$\type(F)=I-i$ for some $i\in I$, and we say that $C$ and $C'$ are
\emph{$i$-adjacent}.  In the canonical example, two distinct chambers
$X_1<X_2<\cdots<X_{n-1}$ and $X'_1<X'_2<\cdots<X'_{n-1}$ are
$i$-adjacent if and only if $X_j=X'_j$ for $j\ne i$.  If we identify
chambers with permutations as above, then $w$ and $w'$ are
$i$-adjacent if and only if the $n$-tuple $(w'(1),w'(2),\dots,w'(n))$
is obtained from $(w(1),w(2),\dots,w(n))$ by interchanging $w(i)$ and
$w(i+1)$.  For example, the chambers labeled 2134 and 2314 in
Figure~\ref{f:braid} (Section~\ref{sub:cell}) are 2-adjacent.

\subsection{Descent sets}

Given two chambers $C,C'$, we define the \emph{descent set} of $C'$
with respect to~$C$, denoted $\des(C,C')$, to be the set of $i\in I$
such that there is a minimal gallery
\[
C=C_0,C_1,\dots,C_l=C'
\]
ending with an $i$-adjacency between $C_{l-1}$ and $C_l$.  (See
section~\ref{sub:gallery} for the definition and basic facts
concerning minimal galleries.)  Equivalently, we have $i\in\des(C,C')$
if and only if $C$ and $C'$ are on opposite sides of the hyperplane
$\supp F$, where $F$ is the face of $C'$ of type $I-i$.  Or, if $C''$
is the chamber $i$-adjacent to~$C'$, then $i\in\des(C,C')$ if and only
if $d(C,C')=d(C,C'')+1$.

If we have chosen a fundamental chamber $C$, then we write $\des(C')$
instead of $\des(C,C')$, and we call it the \emph{descent set}
of~$C'$.  And if $C'$ corresponds to $w\in W$, i.e., if $C'=wC$, then
we also speak of $\des(w)$, the descent set of~$w$.  The terminology
is motivated by the canonical example, where the descent set of a
permutation~$w$ is $\{i:w(i)>w(i+1)\}$; it is a subset
of~$\{1,2,\dots,n-1\}$.  For example, the descent set of 2431 is
$\{2,3\}$; this is consistent with the fact that there are two minimal
galleries from 1234 to 2431 in Figure~\ref{f:braid}, one ending with a
2-adjacency and the other with a 3-adjacency.

Descent sets can be characterized in terms of the semigroup structure
on $\Sigma$:

\begin{proposition} \label{p:descent}
Given chambers $C,C'$ and a face $F\le C'$, we have $FC=C'$ if and
only if $\des(C,C')\subseteq \type(F)$.  Thus $\des(C,C')$ is the type
of the smallest face $F\le C'$ such that $FC=C'$.
\end{proposition}

\begin{proof}
Suppose $FC=C'$.  Given $i\in\des(C,C')$, let $G$ be the face of~$C'$
of type~$I-i$ and let $H=\supp G$; see Figure~\ref{f:descent}.
\begin{figure}[htb]
\setlength{\unitlength}{3947sp}%
\begingroup\makeatletter\ifx\SetFigFont\undefined
\def\x#1#2#3#4#5#6#7\relax{\def\x{#1#2#3#4#5#6}}%
\expandafter\x\fmtname xxxxxx\relax \def\y{splain}%
\ifx\x\y   
\gdef\SetFigFont#1#2#3{%
  \ifnum #1<17\tiny\else \ifnum #1<20\small\else
  \ifnum #1<24\normalsize\else \ifnum #1<29\large\else
  \ifnum #1<34\Large\else \ifnum #1<41\LARGE\else
     \huge\fi\fi\fi\fi\fi\fi
  \csname #3\endcsname}%
\else
\gdef\SetFigFont#1#2#3{\begingroup
  \count@#1\relax \ifnum 25<\count@\count@25\fi
  \def\x{\endgroup\@setsize\SetFigFont{#2pt}}%
  \expandafter\x
    \csname \romannumeral\the\count@ pt\expandafter\endcsname
    \csname @\romannumeral\the\count@ pt\endcsname
  \csname #3\endcsname}%
\fi
\fi\endgroup
\resizebox{!}{2in}{%
\begin{picture}(1887,3237)(1189,-2686)
\thinlines
\put(3001,-976){\circle*{60}}
\put(2101,539){\line( 0,-1){3000}}
\put(2101,-61){\line(-1,-1){900}}
\put(1201,-961){\line( 1,-1){900}}
\put(2101,-1861){\line( 1, 1){900}}
\put(3001,-961){\line(-1, 1){900}}
\put(2426,-1036){\makebox(0,0)[lb]{\smash{\SetFigFont{12}{14.4}{rm}$C'$}}}
\put(1921,-1036){\makebox(0,0)[lb]{\smash{\SetFigFont{12}{14.4}{rm}$G$}}}
\put(3076,-1036){\makebox(0,0)[lb]{\smash{\SetFigFont{12}{14.4}{rm}$i$}}}
\put(2006,-2686){\makebox(0,0)[lb]{\smash{\SetFigFont{12}{14.4}{rm}$H$}}}
\end{picture}}
\caption{}\label{f:descent}
\end{figure}
We know that $C$ and $C'$ are on opposite sides of~$H$, so $F$ must be
strictly on the $C'$-side of~$H$, hence $F\nleq G$ and $i\in\type(F)$.
This proves $\des(C,C')\subseteq \type(F)$.

Conversely, suppose $\des(C,C')\subseteq \type(F)$.  To show $FC=C'$,
it suffices to show that $FC$ and $C'$ are on the same side of every
hyperplane $H=\supp G$, where $G$ is a codimension~1 face of~$C'$; see
\cite[Section I.4B]{\brown}.  This is automatic if $C$ and $C'$ are on
the same side of~$H$, so assume they are not.  Writing $\type(G)=I-i$,
we then have $i\in\des(C,C')$, hence $i\in\type(F)$.  Then $F\nleq G$,
so $F$ is strictly on the $C'$-side of~$H$ and therefore $FC$ is on
the $C'$-side of~$H$.
\end{proof}

\subsection{Descent counts and the $h$-vector} \label{sub:flagh}

In this subsection we fix a fundamental chamber $C$, so that every
chamber $C'$ has a well-defined descent set $\des(C')$.  For
$J\subseteq I$, let $\beta(J)$ be the number of chambers with descent
set~$J$.  This number is independent of the choice of $C$, since the
group of type-preserving automorphisms of~$\Sigma$ is transitive on
the chambers.  It can also be described as the number of $w\in W$ with
descent set $J$.  We will show that the vector $(\beta(J))_{J\subseteq
I}$ coincides with the ``$h$-vector'' defined below; the definition is
modeled on that of the flag $h$-vector for graded posets (see
Section~\ref{sub:flag}).

First we define the \emph{$f$-vector} of $\Sigma$ by setting
$f_J(\Sigma)$ equal to the number of simplices of type~$J$.  The
$h$-vector is then obtained by writing
\begin{equation} \label{e:hvect}
f_J(\Sigma)=\sum_{K\subseteq J} h_K(\Sigma),
\end{equation}
or, equivalently,
\begin{equation}
h_J(\Sigma)=\sum_{K\subseteq J} (-1)^{\abs{J-K}} f_K(\Sigma).
\end{equation}

\begin{proposition} \label{p:beta}
Let $\Sigma$ be the Coxeter complex of a finite reflection group, and
let $I$ be the set of types of vertices.  Then for any $J\subseteq I$,
\[
\beta(J)=h_J(\Sigma).
\]
\end{proposition}

\begin{proof}
Let $\Sigma_J$ be the set of simplices of type $J$.  There is a 1--1
map $\Sigma_J\to\C$, given by $F\mapsto FC$, where $C$ is the
fundamental chamber.  It is 1--1 because we can recover $F$ from $FC$
as the face of type~$J$.  Its image, according to
Proposition~\ref{p:descent}, is the set of chambers with descent set
contained in~$J$.  Hence
\[
f_J(\Sigma)=\sum_{K\subseteq J} \beta(K).
\]
The proposition now follows from \eqref{e:hvect}.
\end{proof}

\begin{remark}
Everything in this and the previous subsection generalizes from finite
Coxeter complexes to finite buildings.
\end{remark}

\subsection{The ring of invariants in the semigroup algebra} \label{sub:inv}

Fix a commutative ring $k$ and consider the semigroup algebra
$k\Sigma$.  This has a natural $W$-action, and the $W$-invariants form
a $k$-algebra $A=(k\Sigma)^W$.  As a $k$-module, $A$ is free with
one basis element for each $W$-orbit in~$\Sigma$, that basis element
being the sum of the simplices in the orbit.  Since orbits correspond
to types of simplices, we get a basis vector
\[
\sigma_J= \sum_{F\in\Sigma_J} F
\]
for each $J\subseteq I$, where, as in the proof of
Proposition~\ref{p:beta}, $\Sigma_J$ is the set of simplices of
type~$J$.  The product of two basis vectors is given by
\[
\sigma_J \sigma_K = \sum_L \alpha_{JKL} \sigma_L,
\]
where $\alpha_{JKL}$ is the number of ways of writing a given simplex
of type $L$ as a product $FG$, where $\type(F)=J$ and $\type(G)=K$.
This number is 0 unless $J\subseteq L$.

There is a second natural basis $(\tau_J)_{J\subseteq I}$ for $A$,
obtained by writing
\[
\sigma_J = \sum_{K\subseteq J} \tau_K,
\]
or, equivalently,
\[
\tau_J = \sum_{K\subseteq J} (-1)^{\abs{J-K}} \sigma_K.
\]
This change of basis is motivated by the study of the $h$-vector
above, and also by considerations in Solomon's paper~\cite{\solomon}.

\subsection{Solomon's descent algebra}

Bidigare \cite{\bid} proved that $A$ is anti-isomorphic to Solomon's
descent algebra, which is a certain subalgebra of the group algebra
$kW$.  We give here a geometric version of his proof.

Recall that the $k$-module $k\C$ spanned by the chambers is an ideal
in $k\Sigma$.  In particular, it is a module over the subring
$A\subseteq k\Sigma$, and the action of $A$ on $k\C$ commutes with the
action of~$W$; we therefore obtain a homomorphism
\[
A\to \End_W(k\C),
\]
the latter being the ring of $kW$-endomorphisms of $k\C$.  

We now choose a fundamental chamber $C\in\C$ and use it to identify
$\C$ with $W$, the correspondence being $wC\leftrightarrow w$.  This is
compatible with left $W$-actions, so $\End_W(k\C)$ gets identified
with the ring of operators on $kW$ that commute with the left action
of~$kW$.  But any such operator $T$ is given by right multiplication
by an element of $kW$, that element being $T(\text{id})$.  So we
obtain, finally, a product-reversing map (i.e., an anti-homomorphism)
\[
\phi\colon A\to kW.
\]
Chasing through the definitions, one sees that $\phi$ is characterized
by
\begin{equation} \label{e:phi}
\phi(a)C=aC
\end{equation}
for $a\in A$.  Here $C$ is the fundamental chamber, the product on the
left is given by the action of $W$ on $\C$, and the product on the
right takes place in the semigroup algebra~$k\Sigma$.

Recall now that the choice of fundamental chamber determines a special
set of generators $S$ for $W$ (the ``simple reflections'', or
``Coxeter generators''), consisting of the reflections with respect to
the supports of the codimension~1 faces of~$C$.  Using this set of
generators, Solomon~\cite{\solomon} defined a subalgebra of $kW$,
which has come to be known as the \emph{descent algebra}.

\begin{theorem}[Bidigare] \label{t:descent}
Let $W$ be a finite reflection group with Coxeter complex $\Sigma$,
and let $A$ be the invariant subalgebra $(k\Sigma)^W$ of the semigroup
algebra $k\Sigma$, where $k$ is an arbitrary commutative ring.  Then
$A$ is anti-isomorphic to Solomon's descent algebra.
\end{theorem}

\begin{proof}
Let $(\sigma_J)$ and $(\tau_J)$ be the $k$-bases of $A$ introduced in
Section~\ref{sub:inv}.  We have
\[
\sigma_J C = \sum_{F\in\Sigma_J} FC.
\]
As we noted in the proof of Proposition \ref{p:beta}, the chambers
$FC$ that occur in the sum are those with descent set contained in
$J$.  Under our bijection between $\C$ and $W$, a given chamber $FC$
corresponds to the element $w\in W$ such that $FC=wC$.  So we can
write
\[
\sigma_J C = \sum_{w\in U_J} wC,
\]
where $U_J = \{w\in W:\des(w)\subseteq J\}$.  Our characterization
\eqref{e:phi} of $\phi$ therefore yields
\[
\phi(\sigma_J)=u_J:=\sum_{w\in U_J} w.
\]
Let $Z_J=\{w\in W:\des(w)=J\}$ and let $z_J=\sum_{w\in Z_J} w$.  Then
$u_J=\sum_{K\subseteq J} z_K$, so
\[
\sum_{K\subseteq J} \phi(\tau_K) = \sum_{K\subseteq J} z_K
\]
for all $J$, hence $\phi(\tau_J)=z_J$.  Since the $z_J$ are clearly
linearly independent, it follows that $\phi$ is injective and hence
gives an anti-isomorphism of $A$ onto a subalgebra of~$kW$.

It remains to show that $\phi(A)$ is the descent algebra.  For each
$i\in I$, let $s_i$ be the reflection with respect to the face of~$C$
of type $I-i$.  Then our set of generators of~$W$ is $S=\{s_i:i\in
I\}$.  Moreover, for any $J\subseteq I$ the stabilizer of the face
of~$C$ of type~$J$ is the subgroup $W_{I-J}$ generated by $\{s_i:i\in
I-J\}$; for example, the face of type $I-i$ has stabilizer of order~2,
generated by $s_i$.  Finally, our definition of descent sets has the
following translation: $i\in\des(w)$ if and only if
$\ell(ws_i)<\ell(w)$, where $\ell$ is the length function on~$W$ with
respect to the generating set~$S$.  Using these remarks, the reader
can easily check that our basis vector $u_J$ for $\phi(A)$ coincides
with Solomon's $x_T$, where $T=\{s_i : i\in I-J\}\subseteq S$.  Hence
$\phi(A)$ is equal to the descent algebra.
\end{proof}

\begin{remark}
There has been interest recently in giving explicit formulas for
orthogonal families of idempotents in the descent algebra that lift
the primitive idempotents of the algebra mod its radical.  See, for
example, Bergeron et al \cite{\bergeron} and earlier references cited
there.  The results of the present paper provide further formulas of
this type; it suffices to take a generic element $p=\sum_{F\in\Sigma}
p_F F \in (k\Sigma)^W$, find the primitive idempotents in $k[w]$ by
the results of Section~\ref{s:semisimple}, and apply~$\phi$.
\end{remark}

\subsection{The descent algebra and random walks} \label{sub:descentwalk}

Assume, for the moment, that we have \emph{not} chosen a fundamental
chamber.  Given a $W$-invariant probability distribution
$p=\sum_{F\in\Sigma} p_F F$, we get a $W$-invariant random walk on
$\C$ with transition matrix
\[
P(C,D) = \sum_{FC=D} p_F
\]
for $C,D\in\C$.  (``$W$-invariant'' means $P(wC,wD)=P(C,D)$ for
$C,D\in\C$, $w\in W$.)  If we now choose a fundamental chamber $C$ and
identify $\C$ with $W$, we get a left-invariant Markov chain on~$W$
whose transition matrix satisfies
\[
P(\text{id},w) = \mu_w:=\sum_{FC=wC} p_F.
\]
Left invariance implies that this is a \emph{right} random walk on the
group $W$:  At each step, we choose $w$ with probability $\mu_w$ and
right-multiply by~$w$.  Note that the definition of the probability
distribution $\{\mu_w\}$ on~$W$ can be written as $\mu C=pC$, where
$\mu=\sum_{w\in W} \mu_w w \in\R W$; as in~\eqref{e:phi}, the product
on the left is given by the action of~$W$ and the product on the right
is in~$\R\Sigma$.  Hence \eqref{e:phi} implies that $\mu=\phi(p)$.  We
have therefore proved:

\begin{theorem} \label{t:descentwalk}
Let $W$ be a finite reflection group with Coxeter complex $\Sigma$,
and let $p$ be a $W$-invariant probability distribution on~$\Sigma$.
Choose a fundamental chamber and use it to identify $\C$ with~$W$.
Then the hyperplane chamber walk on $\C$ driven by~$p$ corresponds to
the right random walk on $W$ driven by $\mu$, where $\mu\in\R W$ is
the image of $p$ under the isomorphism of Theorem~\ref{t:descent}
between $(\R\Sigma)^W$ and the descent subalgebra of $\R W$.
Consequently, the algebra $\R[\mu]\subseteq \R W$ generated by $\mu$
is a split semisimple commutative subalgebra of the descent algebra.
\end{theorem}

Suppose, for example, that $p$ is uniform on simplices of type $J$ for
some fixed~$J$.  Then the proof of Theorem~\ref{t:descent} shows that
$\mu=\phi(p)$ is uniform on the $w\in W$ with $\des(w)\subseteq J$.
This explains some of the observations in \cite{\riffle,\dfp}.
Returning to our canonical example with $W=S_n$, let $J=\{1\}$; thus
$p$ is uniform on the vertices of type~1, i.e., the singleton subsets
of~$[n]$.  The corresponding hyperplane chamber walk is the Tsetlin
library with uniform weights.  Viewing this as a walk on the
permutation group $S_n$, it is the right random walk driven by the
uniform distribution $\mu$ on the permutations $w$ with
$\des(w)\subseteq\{1\}$.  There are $n$ such, with
$(w(1),\dots,w(n))=(i,1,2,\dots,i-1,i+1,\dots,n)$, $i=1,\dots,n$.

Continuing with this example, we can use formula~\eqref{e:tsetlinu}
(Section~\ref{sub:tsetlinu}) to get formulas for the primitive
idempotents in $\R[\mu]$.  In fact, the probability
distribution~$\bar{w}\in\R\bar{S}$ of Section~\ref{sub:tsetlinu} is
the same as what we are now calling~$p$.  (Recall from
Section~\ref{sub:fnbar} that $\bar{S}\subset\Sigma$, so this assertion
makes sense.)  And the element~$\bar{\sigma}_l$ of
Section~\ref{sub:tsetlinu} is the same as the element
$\sigma_{\{1,\dots,l\}}\in(\R \Sigma)^W$ if $l<n$, while
$\bar{\sigma_n}=\bar{\sigma}_{n-1}$.  Combining the isomorphism
$\R[\bar{w}] = \R[p]\isoto \R[\mu]$ with equation~\eqref{e:tsetlinu},
we now obtain the following result:  Define $v_l\in\R W$ for $0\le
l\le n$ by
\[
v_l=\sum_{\des(w)\subseteq\{1,\dots,l\}} w \qquad\text{if $l<n$},
\]
and
\[
v_n=v_{n-1}.
\]
Let
\[
E_i= \sum_{l=i}^n (-1)^{l-i} \binom{l}{i}\frac{v_l}{l!}.
\]
Then~$E_0,\dots,E_{n-2},E_n$, are the primitive idempotents
in~$\R[\mu]$.  These formulas are the same as those of
\cite[Theorem~4.2]{\dfp}.

\appendix
\section{The hyperplane face semigroup} \label{app:hyperplane}

More details concerning the material reviewed here can be found in
\cite{\bhr,\bbd,\red,\brown,\bd,\ot,\zie}.  Throughout this section
$\A=\{H_i\}_{i\in I}$ denotes a finite set of affine hyperplanes in
$V=\R^n$.  Let $H_i^+$ and $H_i^-$ be the two open halfspaces
determined by $H_i$; the choice of which one to call $H_i^+$ is
arbitrary but fixed.

\subsection{Faces and chambers}

The hyperplanes $H_i$ induce a partition of $V$ into convex sets
called \emph{faces} (or \emph{relatively open faces}).  These are the
nonempty sets $F \subseteq V$ of the form
\[
F = \bigcap_{i\in I} H_i^{\sigma_i},
\]
where $\sigma_i \in \{+,-,0\}$ and $H^0_i = H_i$.  Equivalently, if we
choose for each $i$ an affine function $f_i\colon V \rightarrow \R$
such that $H_i$ is defined by $f_i = 0$, then a face is a nonempty set
defined by equalities and inequalities of the form $f_i > 0$, $f_i <
0$, or $f_i = 0$, one for each $i \in I$.  The sequence $\sigma =
(\sigma_i)_{i\in I}$ that encodes the definition of $F$ is called the
\emph{sign sequence} of $F$ and is denoted $\sigma(F)$.

The faces such that $\sigma_i \ne 0$ for all $i$ are called
\emph{chambers}.  They are convex open sets that partition the
complement $V - \bigcup_{i\in I} H_i$.  In general, a face $F$ is
open relative to its \emph{support}, which is defined to be the affine
subspace
\[
\supp F = \bigcap_{\sigma_i (F) = 0} H_i.
\]
Since $F$ is open in $\supp F$, we can also describe $\supp F$ as the
affine span of~$F$.

\subsection{The face relation} 

The \emph{face poset} of $\A$ is the set $\F$ of faces, ordered as
follows:  $F \le G$ if for each $i \in I$ either $\sigma_i(F) = 0$ or
$\sigma_i(F) = \sigma_i(G)$.  In other words, the description of $F$
by linear equalities and inequalities is obtained from that of $G$ by
changing zero or more inequalities to equalities.

\subsection{Product}
The set $\F$ of faces is also a semigroup.  Given $F,G \in \F$, their
\emph{product} $FG$ is the face with sign sequence
\[
\sigma_i(FG) = 
\begin{cases}
  \sigma_i(F)  &\text{if $\sigma_i(F) \ne 0$} \\
         \sigma_i(G) &\text{if $\sigma_i(F) = 0$.}
\end{cases}
\]
This has a geometric interpretation:  If we move on a straight line
from a point of~$F$ toward a point of $G$, then $FG$ is the face we
are in after moving a small positive distance.  Notice that the face
relation can be described in terms of the product:  One has
\begin{equation} \label{e:poset}
F\le G \iff FG=G.
\end{equation}

\subsection{The semilattice of flats}

A second poset associated with the arrangement $\A$ is the
\emph{semilattice of flats}, also called the \emph{intersection
semilattice}, which we denote by~$\L$.  It consists of all nonempty
affine subspaces $X \subseteq V$ of the form $X = \bigcap_{H\in\A'}
H$, where $\A' \subseteq \A$ is an arbitrary subset (possibly empty).
We order $\L$ by inclusion.  [Warning: Many authors order $\L$ by
reverse inclusion.]  Notice that any two elements $X,Y$ have a
least upper bound $X\vee Y$ in~$\L$, which is the intersection of
all hyperplanes $H\in\A$ containing both $X$ and~$Y$; hence $\L$
is an \emph{upper semilattice} (poset with least upper bounds).  It is
a lattice if the arrangement $\A$ is \emph{central}, i.e., if
$\bigcap_{H\in\A} H \ne\emptyset$.  Indeed, this intersection is then
the smallest element of~$\L$, and a finite upper semilattice with a
smallest element is a lattice \cite[Section 3.3]{\stanley}.  The
support map gives a surjection
\[
\supp\colon\F\onto\L,
\]
which preserves order and also behaves nicely with respect to the
semigroup structure.  Namely, we have
\begin{equation}
\supp(FG) = \supp F \vee \supp G
\end{equation}
and
\begin{equation}
FG = F \iff \supp G \le \supp F.
\end{equation}

\subsection{Example: The braid arrangement} \label{sub:braid}

The \emph{braid arrangement} in $\R^n$ consists of the $\binom{n}{2}$
hyperplanes $H_{ij}$ defined by $x_i=x_j$, where $1\le i<j\le n$.
Each chamber is determined by an ordering of the coordinates, so it
corresponds to a permutation.  When $n=4$, for example, one of the 24
chambers is the region defined by $x_2>x_3>x_1>x_4$, corresponding to
the permutation 2314.  The faces of a chamber~$C$ are obtained by
changing to equalities some of the inequalities defining~$C$.  For
example, the chamber $x_2>x_3>x_1>x_4$ has a face given by
$x_2>x_3>x_1=x_4$, which is also a face of the chamber
$x_2>x_3>x_4>x_1$.

It is useful to encode the system of equalities and inequalities
defining a face~$F$ by an ordered partition $(B_1,\dots,B_k)$ of
$[n]=\{1,\dots,n\}$.  Here $B_1,\dots,B_k$ are disjoint nonempty sets
whose union is $[n]$, and their order counts.  For example, the face
$x_2>x_3>x_1=x_4$ corresponds to the 3-block ordered partition
$(\{2\},\{3\},\{1,4\})$, and the face $x_2>x_1=x_3=x_4$ corresponds to
the 2-block ordered partition $(\{2\}, \{1,3,4\})$.

Thus the face semigroup of the braid arrangement can be viewed as the
set $\B$ of ordered partitions, with a product that one can easily
work out.  We have recorded this product in Section~\ref{sub:braidsemi},
where one can also find a description of the face relation, the
intersection lattice, and the support map.  See also
Section~\ref{sub:reflection}, where the braid arrangement appears as
the canonical example of a reflection arrangement.

\subsection{Spherical representation} \label{sub:cell}

Suppose now that $\A$ is a \emph{central} arrangement, i.e., that the
hyperplanes have a nonempty intersection.  We may assume that this
intersection contains the origin.  Suppose further that $\bigcap_{i\in
I} H_i = \{0\}$, in which case $\A$ is said to be \emph{essential}.
(There is no loss of generality in making this assumption; for if it
fails, then we can replace $V$ by the quotient space $V/\bigcap_i H_i$.)
The hyperplanes then induce a cell-decomposition of the unit sphere,
the cells being the intersections with the sphere of the faces
$F\in\F$.  Thus $\F$, as a poset, can be identified with the poset of
cells of a regular cell-complex $\Sigma$, homeomorphic to a sphere.
Note that the face $F=\{0\}$, which is the identity of the
semigroup~$\F$, is not visible in the spherical picture; it
corresponds to the empty cell.  The cell-complex $\Sigma$ plays a
crucial role in~\cite{\bd}, to which we refer for more details.

\begin{figure}[htb]
\begin{center}
\resizebox{!}{2.5in}{\includegraphics{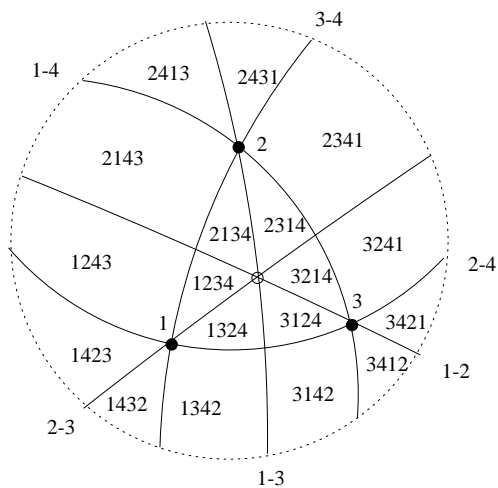}}
\end{center}
\caption{The braid arrangement when $n=4$.} \label{f:braid}
\end{figure}

The braid arrangement provides a simple example.  It is not essential,
because the hyperplanes $H_{ij}$ intersect in the line $L$ defined by
$x_1=\cdots=x_n$.  We can therefore view the braid arrangement as an
arrangement in the $(n-1)$-dimensional quotient space $\R^n/L$.  When
$n=4$, we obtain an arrangement of six planes in~$\R^3$, whose
spherical picture is shown in Figure~\ref{f:braid}.
The plane corresponding to $H_{ij}$ cuts the sphere in the great
circle labeled i-j.  Each chamber of the arrangement is a
simplicial cone, which intersects the sphere in a triangle labeled
with the associated permutation.  Figure \ref{f:braid} has been
reproduced from~\cite{\bbd}, where one can find further discussion and
more examples.

\subsection{Galleries and convex sets} \label{sub:gallery}

We return to an arbitrary arrangement $\A$.  Two chambers $C,C'\in\C$
are said to be \emph{adjacent} if they have a common codimension~1
face.  A \emph{gallery} is a sequence of chambers $C_0,C_1,\dots,C_l$
such that $C_{j-1}$ and $C_j$ are adjacent for each $j=1,2,\dots,l$.
Given $C,C'\in\C$, the minimal length $l$ of a gallery from $C$
to~$C'$ is the \emph{distance} between $C$ and~$C'$, denoted
$d(C,C')$; and any gallery from $C$ to~$C'$ of minimal length
$d(C,C')$ is called a \emph{minimal gallery}.  The distance $d(C,C')$
can also be characterized as the number of hyperplanes in~$\A$
separating $C$ from~$C'$; in fact, every minimal gallery from $C$
to~$C'$ crosses each of these hyperplanes exactly once (see
\cite[Section~I.4E]{\brown}).

Let $\D\subseteq\C$ be a nonempty set of chambers.  We say that $\D$
is \emph{convex} if it satisfies the equivalent conditions of the
following result:

\begin{proposition} \label{p:convex}
The following conditions on a nonempty set $\D\subseteq\C$ are
equivalent:
\begin{enumerate}
\item[(i)] 
For any $C,C'\in\D$, every minimal gallery from $C$ to~$C'$
is contained in~$\D$.
\item[(ii)] 
$\D$ is the set of chambers in an intersection of some of
the halfspaces determined by~$\A$.
\end{enumerate}
\end{proposition}

In terms of sign sequences, condition (ii) says that there is a subset
$J\subseteq I$ and a set of signs $\sigma_i\in\{+,-\}$, $i\in J$, such
that
\[
\D=\{C\in\C : \sigma_i(C)=\sigma_i \text{ for all } i\in J\}.
\]

Proposition~\ref{p:convex}, which is stated as in exercise in
\cite[Section~I.4E]{\brown}, is essentially due to Tits
\cite[Theorem~2.19]{\tits}.  See also \cite[Proposition~4.2.6]{\red}
for a proof in the context of oriented matroids.  For the convenience
of the reader, here is the latter proof specialized to hyperplane
arrangements:

\begin{proof}
A minimal gallery from $C$ to $C'$ crosses only the hyperplanes that
separate $C$ from~$C'$.  This shows that (ii) implies~(i).  For the
converse, it suffices to show that if (i) holds and $C$ is a chamber
not in~$\D$, then there is a hyperplane $H\in\A$ separating $C$
from~$\D$.  Choose a gallery $D,C_1,C_2,\dots,C_l=C$ of minimal
length, starting in~$\D$ and ending at~$C$.  By minimality, we have
$C_1\notin\D$.  Let $H$ be the (unique) hyperplane in~$\A$ separating
$D$ from~$C_1$.  Then $H$ also separates $D$ from~$C$.  For any
$D'\in\D$, we have $d(D,D')=d(C_1,D')\pm1$, where the sign depends on
which $H$-halfspace contains~$D'$.  The sign cannot be $+$, because
then we could construct a minimal gallery from $D$ to~$D'$ passing
through~$C_1$, contradicting~(i).  So $d(D,D')=d(C_1,D')-1$, which
means that $D$ and~$D'$ are on the same side of~$H$.  Thus $H$
separates $\D$ from~$C$, as required.
\end{proof}

\section{Left-regular bands: Foundations} \label{app:semi}

In this appendix $S$ is an arbitrary semigroup, not necessarily
finite, not necessarily having an identity.  Motivated by the theory
of hyperplane face semigroups, we wish to isolate the conditions
on~$S$ under which we can define analogues of the face relation,
chambers, the semilattice of flats, etc.

\subsection{Partial order}
Given $x,y\in S$, we set $x\le y$ if $xy=y$.  This relation is
transitive for any semigroup (see the first paragraph
of~Section~\ref{sub:poset}).  It is reflexive if and only if every
element of $S$ is idempotent, in which case $S$ is called an
\emph{idempotent semigroup} or a \emph{band}.  Antisymmetry, however,
imposes a much stronger condition on~$S$:

\begin{proposition}
The relation defined above is a partial order if and only if $S$ is an
idempotent semigroup satisfying
\begin{equation} \label{e:regular1}
xyx=xy
\end{equation}
for all $x,y\in S$.
\end{proposition}

In other words, the relation makes $S$ a poset if and only if $S$ is a
LRB as defined in Section~\ref{sub:semi}.

\begin{proof}
For the ``if'' part see the beginning of Section \ref{sub:poset}.  To
prove the converse, we may assume that $S$ is an idempotent semigroup
for which the relation is antisymmetric, and we must
prove~\eqref{e:regular1}.  Note that $(xy)(xyx)=(xy)^2x=xyx$, so
$xy\le xyx$.  On the other hand, $(xyx)(xy)=xyxy=(xy)^2=xy$, so
$xyx\le xy$.  Thus antisymmetry implies that $xyx=xy$, as required.
\end{proof}

\subsection{The associated semilattice}

We now show how to construct, for any idempotent semigroup
satisfying~\eqref{e:regular1}, an analogue of the intersection
semilattice of a hyperplane arrangement.  In particular, this shows
that the definition of LRB given in Section~\ref{s:intro} is
equivalent to the one given in Section~\ref{s:lrb} and used throughout
this paper.

\begin{proposition}
Let $S$ be an idempotent semigroup satisfying \eqref{e:regular1}.
Then there is a semilattice~$L$ that admits an order-preserving
surjection $\supp\colon S\onto L$ such that
\begin{equation} \label{e:homom1}
\supp xy = \supp x \vee \supp y
\end{equation}
for all $x,y\in S$ and
\begin{equation} \label{e:delete1}
xy=x \iff \supp y\le\supp x.
\end{equation}
\end{proposition}

\begin{proof}
The construction of $L$ is forced on us by \eqref{e:delete1}:  Define
a relation $\preceq$ on $S$ by $y\preceq x \iff xy=x$.  This is
transitive and reflexive, but not necessarily antisymmetric.  We
therefore obtain a poset $L$ by identifying $x$ and~$y$ if $x\preceq
y$ and $y\preceq x$.  If we denote by $\supp\colon S\onto L$ the
quotient map, then \eqref{e:delete1} holds by definition.  To see that
$\supp$ is order-preserving, suppose that $x\le y$, i.e., $xy=y$.
Multiplying on the right by~$x$ and using \eqref{e:regular1}, we
conclude that $xy=yx$; hence $yx=y$ and $x\preceq y$, i.e., $\supp
x\le \supp y$.  It remains to show that $\supp xy$ is the least upper
bound of $\supp x$ and $\supp y$ in~$L$.  It is an upper bound because
the equations $xyx=xy$ and $xyy=xy$ show that $xy\succeq x$ and
$xy\succeq y$.  And it is the least upper bound, because if $z\succeq
x$ and $z\succeq y$, then $zx=z$ and $zy=z$, whence
$z(xy)=(zx)y=zy=z$, so that $z\succeq xy$.
\end{proof}

If $S$ has an identity e, then $L$ has a smallest element
$\zerohat=\supp e$.  If, in addition, $L$ is finite, then it is is a
lattice \cite[Section 3.3]{\stanley}.

\subsection{Chambers}

We close this appendix by giving several characterizations of the
chambers.  Let $S$ be a LRB whose semilattice $L$ has a largest
element~$\onehat$.  This is automatic if $S$ is finite.  As in
Section~\ref{sub:semi}, we call an element $c\in S$ a \emph{chamber}
if $\supp c=\onehat$.

\begin{proposition}
The following conditions on an element $c\in S$ are equivalent:
\begin{enumerate}
\item[(i)]
$c$ is a chamber.
\item[(ii)]
$cx=c$ for all $x\in S$.
\item[(iii)]
$c$ is maximal in the poset $S$.
\end{enumerate}
\end{proposition}

\begin{proof}
We have $\supp c=\onehat \iff \supp c\ge\supp x$ for all $x\in S$.  In
view of~\eqref{e:delete1}, this holds if and only if $cx=x$,
so (i) and~(ii) are equivalent.  If (ii) holds then $c$ is maximal,
because $c\le x \implies cx=x \implies c=x$.  For the converse, note
that $c\le cx$ for all $x,c\in S$; so if $c$ is maximal then (ii)
holds.
\end{proof}

The set $C$ of chambers is a 2-sided ideal in $S$.  Indeed, (ii) shows
that it is a right ideal, and it is a left ideal because if $\supp
c=\onehat$ then $\supp xc =\onehat$ by~\eqref{e:homom1}.  One can
check that $C$ is the \emph{kernel} of the semigroup~$S$, i.e., the
(unique) minimal 2-sided ideal.

\section{Generalized derangement numbers} \label{app:derangement}

In this appendix we associate to any finite poset~$L$ with
$\zerohat,\onehat$ a \emph{derangement number}~$d(L)\ge0$.  If $L$ is
the Boolean lattice of rank~$n$, then $d(L)$ is the ordinary
derangement number ~$d_n$ (number of fixed-point-free permutations of
an $n$-set).  If $L$ is the lattice of subspaces of an $n$-dimensional
vector space over $\fq$, then $d(L)$ is the $q$-analogue of~$d_n$
studied by Wachs~\cite{\wachs}.  If $L$ is the lattice of contractions
of a graph, then $d(L)$ is some (new) graph invariant.

We are mainly interested in the case where $L$ is a \emph{geometric
lattice}, i.e., the lattice of flats of a matroid.  In this case, the
derangement numbers of the intervals $[X,\onehat]$ give the
multiplicities of the eigenvalues for the random walk on the maximal
chains of~$L$ constructed in Section~\ref{sub:matroidsemi}.  But since
the derangement numbers may be of independent interest, we will keep
this appendix logically independent of the theory of random walks; the
latter will be mentioned only for motivation.

\subsection{Definition} \label{sub:defn1}

Let $L$ be a finite poset with smallest element $\zerohat$ and largest
element~$\onehat$.  We associate to $L$ an integer $d(L)$, called the
\emph{derangement number} of~$L$.  It is defined inductively by the
equation
\begin{equation} \label{e:recurrence}
\sum_{X\in L} d([X,\onehat]) = f(L),
\end{equation}
where $f(L)$ is the number of maximal chains in~$L$.  If $L=0$
(the one-element poset, with $\zerohat=\onehat$), this gives $d(L)=1$.
Otherwise, it gives a recurrence that can be solved for
$d(L)=d([\zerohat,\onehat])$; thus
\begin{equation} \label{e:recurrence1}
d(L) = f(L) - \sum_{X>\zerohat} d([X,\onehat]).
\end{equation}
Note that $d(L)=0$ if $L$ is the two-element poset
$\{\zerohat,\onehat\}$.  More generally, $d(L)=0$ if $L$ has exactly
one atom, where an \emph{atom} is a minimal element of $L-\zerohat$.
Indeed, let $X_0$ be the atom and let $L_0=[X_0,\onehat]$.  Then
$f(L)=f(L_0)$, so \eqref{e:recurrence1} becomes
\[
d(L) = f(L_0) -\sum_{X\in L_0} d([X,\onehat]),
\]
and the right side is 0 by \eqref{e:recurrence} applied to $L_0$.

If we apply the definition \eqref{e:recurrence} to each interval
$[Y,\onehat]$, we get
\begin{equation} \label{e:recurrencey}
f([Y,\onehat])=\sum_{X\ge Y}d([X,\onehat]).
\end{equation}
In case $L$ is a geometric lattice, this system of equations for the
numbers $d([X,\onehat])$ is the same as the system of
equations~\eqref{e:mult3} in Section~\ref{sub:matroidsemi} for the
multiplicities $m_X$; this proves our assertion that
$d([X,\onehat])=m_X$.  And this interpretation of $d([X,\onehat])$
also provides an easy way to remember the
definition~\eqref{e:recurrence}, which says that the sum of the
multiplicities equals the size of the state space for the random walk.

We can solve \eqref{e:recurrencey} by M\"{o}bius inversion to get
\[
d([Y,\onehat])=\sum_{X\ge Y}\mu(Y,X)f([X,\onehat]).
\]
Setting $Y=\zerohat$, we get an explicit formula for $d(L)$:
\begin{equation} \label{e:mobius}
d(L)=\sum_{X\in L} \mu(\zerohat,X)f([X,\onehat]).
\end{equation}
It is useful to have a slight variant of this:
\begin{equation} \label{e:recurrence2}
d(L)= \mu(\zerohat,\onehat)+\sum_{X\in M} d([\zerohat,X]),
\end{equation}
where $M$ is the set of maximal elements of $L-\onehat$.  This is
proved by writing
\[
f([X,\onehat])=\sum_{\substack{Y\in M\\Y\ge X}} f([X,Y])
\]
for $X<\onehat$, and then rearranging the sum in \eqref{e:mobius}.

It is not clear from what we have done so far that $d(L)\ge0$, though
we know this is true if $L$ is geometric, since it is the multiplicity
$m_{\zerohat}$.  An independent proof of this, valid for any $L$, is
obtained by giving yet another recursive formula for~$d(L)$, which
involves no signs.

\begin{proposition} \label{p:defn}
If $L=0$ then $d(L)=1$.  Otherwise,
\begin{equation} \label{e:defn}
d(L)=\sum_{X<\onehat} (c(X)-1)d([\zerohat,X]),
\end{equation}
where $c(X)$ is the number of covers of~$X$.
\end{proposition}

(Recall that $Y$ \emph{covers} $X$, written $X\lessdot Y$, if $X<Y$
and there is no $Z$ with $X<Z<Y$.)

\begin{corollary}
$d(L)\ge0$, with equality if and only if $L$ has exactly one atom.
\end{corollary}

\begin{proof}[Proof of the corollary]
The inequality is immediate by induction on the size of $L$.  We have
already observed that equality holds if $L$ has exactly one atom.  If
$L$ has no atoms, then $L=0$ and $d(L)=1>0$.  If $L$ has more than one
atom, then consideration of the term $X=\zerohat$ in~\eqref{e:defn}
shows that $d(L)>0$.
\end{proof}

\begin{proof}[Proof of Proposition \ref{p:defn}]
Let us temporarily take the statement of the proposition as a new
definition of $d(L)$.  It then suffices to show that, with this
definition, equation~\eqref{e:recurrence} holds.  We may assume that
$L\ne0$ and that \eqref{e:recurrence} holds for smaller posets.  Then
\begin{align*}
\sum_{X\in L} d([X,\onehat]) &= 1+ \sum_{X<\onehat}d([X,\onehat])\\
&= 1+\sum_{X<\onehat} \sum_{X\le Y<\onehat}(c(Y)-1)d([X,Y]) &
&\text{by \eqref{e:defn}}\\
&=1+\sum_{Y<\onehat}(c(Y)-1)\sum_{X\le Y}d([X,Y])&&\\
&=1+\sum_{Y<\onehat}(c(Y)-1)f([\zerohat,Y])&&\text{by induction}\\
&=1+\sum_{Y\in L}(c(Y)-1)f([\zerohat,Y]) + f(L).
\end{align*}
So we are done if we can show $1+\sum_{Y\in L}
(c(Y)-1)f([\zerohat,Y])=0$, i.e.,
\[
1+\sum_{Y\in L}c(Y)f([\zerohat,Y]) = \sum_{Y\in L}f([\zerohat,Y]).
\]
The sum on the right counts all chains $\zerohat=X_0\lessdot X_1
\lessdot \cdots \lessdot X_m$ in~$L$, where $m\ge0$.  The sum on the
left counts all such chains of length $m>0$.  Adding 1 counts the
chain of length~0, so the equation holds.
\end{proof}
 
\subsection{Examples} \label{sub:examples}
\setcounter{example}{0}

\begin{example}[Ordinary derangement numbers]
Let $L$ be the Boolean lattice of subsets of an $n$-set.  Writing
$d(L)=d_n$, the recurrence~\eqref{e:recurrence} becomes
\[
\sum_{i=0}^n \binom{n}{i}d_i=n!,
\]
which is a standard recurrence for the ordinary derangement numbers.
(It is obtained by counting permutations according to the number of
elements they move.)  Formulas \eqref{e:mobius}
and~\eqref{e:recurrence2} are the well-known results
\[
d_n=\sum_{i=0}^n (-1)^i\binom{n}{i}(n-i)! =
n!\sum_{i=0}^n\frac{(-1)^i}{i!}
\]
and
\begin{equation} \label{e:recurrence3}
d_n= nd_{n-1}+(-1)^n;
\end{equation}
see~\cite[Section 2.2]{\stanley}.  Finally, Proposition \ref{p:defn} reads
\begin{align*}
d_0&=1\\
d_n&=\sum_{i=0}^{n-1}\binom{n}{i}(n-i-1)d_i \qquad(n>0),
\end{align*}
which may be new.
\end{example}

\begin{example}[$q$-analogue] \label{ex:qder}
Let $L$ be the lattice of subspaces of an $n$-dimensional vector space
over $\fq$.  Writing $d(L)=d_n$ [$=d_n(q)$], the
recurrence~\eqref{e:recurrence} becomes
\[
\sum_{i=0}^n \qbinom{n}{i}d_i=[n]!,
\]
which characterizes the $q$-derangement numbers of Wachs
\cite[p.~277]{\wachs}.  Here $\qbinom{n}{i}$ and $[n]!$ are the
$q$-analogues of $\binom{n}{i}$ and $n!$, respectively.  The inverted
form of this as in~\eqref{e:mobius} is
\[
d_n=\sum_{i=0}^n (-1)^i\qbinom{n}{i}[n-i]!q^{\binom{i}{2}}=
[n]!\sum_{i=0}^n\frac{(-1)^i}{[i]!}q^{\binom{i}{2}};
\]
see~\cite[Theorem 4]{\wachs}.  Finally, Proposition \ref{p:defn} reads
\begin{align*}
d_0&=1\\
d_n&=\sum_{i=0}^{n-1}\qbinom{n}{i}([n-i]-1)d_i \qquad(n>0),
\end{align*}
where $[n-i]$ is the $q$-analogue of $n-i$.
\end{example}

\begin{example}[A graph invariant]
Let $L=L(G)$ be the lattice of contractions of a simple finite graph
$G$, as discussed in Section~\ref{sub:graph}.  Set
$f(G)=f(L(G))$ and $d(G)=d(L(G))$.  Thus $f(G)$ is the number of
collapsing sequences of~$G$ and $d(G)$ is some new invariant of~$G$,
defined by
\[
\sum_{\bar{G}} d(\bar{G}) = f(G),
\]
where the sum is taken over all collapsings $\bar{G}=G/\Pi$.  The
inverted form is
\[
d(G)=\sum_{\Pi\in L(G)} \mu(\zerohat,\Pi)f(G/\Pi).
\]
The numbers $\mu(\zerohat,\Pi)$ that occur here are familiar from
Rota's formula for the chromatic polynomial of~$G$~\cite[Chapter
3, Exercise 44]{\stanley}:  One has
\[
\chi_G(x)=\sum_{\Pi\in L(G)}\mu(\zerohat,\Pi)x^{\abs{\Pi}}.
\]
Finally, Proposition \ref{p:defn} gives
\begin{align*}
d(G)&=1 &&\text{if $G$ is discrete}\\
d(G)&=\sum_{\Pi<\onehat}(e(G/\Pi)-1) d(G_\Pi) &&\text{otherwise},
\end{align*}
where $e(\enspace)$ denotes the number of edges of a graph and
$G_\Pi\subset G$ is the union of the subgraphs induced by the blocks
of~$\Pi$.
\end{example}

\subsection{Connection with the flag $h$-vector} \label{sub:flag}

The result of this subsection (Proposition \ref{p:stanley}) is due
to Richard Stanley and is included with his permission.

The ordinary derangement number $d_n$ has the following
interpretation, due to D\'esarm\'enien~\cite{\desar} (see also
\cite{\desarw} and further references cited there):  Call a
permutation $\pi\in S_n$ a \emph{des\-arrangement} if the maximal
initial descending sequence $\pi(1)>\pi(2)>\cdots>\pi(l)$ has even
length~$l$; then $d_n$ is the number of desarrangements.
D\'esarm\'enien gave a bijective proof of this assertion and used it
to give a combinatorial proof of the recurrence \eqref{e:recurrence3}.
One can also reverse the process and deduce D\'esarm\'enien's result
from~\eqref{e:recurrence3}, by induction on~$n$.

The result can be phrased in terms of descent sets.  Recall that $\pi$
is said to have a \emph{descent} at~$i$ if $\pi(i)>\pi(i+1)$, where
$1\le i\le n-1$.  For $J\subseteq[n-1]=\{1,\dots,n-1\}$, let
$\beta(J)$ be the number of permutations in $S_n$ with descent set
$J$.  Let $\J$ be the family of sets $J$ such that the first integer
$l\ge1$ not in~$J$ is even.  Then D\'esarm\'enien's interpretation of
$d_n$ is
\begin{equation} \label{e:descent}
d_n=\sum_{\substack{J\subseteq[n-1]\\ J\in\J}} \beta(J).
\end{equation}
We wish to generalize this.  The role of the descent numbers
$\beta(J)$ is played by the components of the flag $h$-vector.
We briefly recall the definition of the latter; for more information,
see \cite[Section III.4]{\stanleyflag}, \cite[Sections 3.12 and
3.8]{\stanley}, or~\cite{\bl}.

Let $L$ be a graded poset with $\zerohat,\onehat$; thus all maximal
chains have the same length $n$, called the \emph{rank} of~$L$.  For
$J\subseteq[n-1]$, let $f_J(L)$ be the number of flags in~$L$ of
type~$J$, where the \emph{type} of a flag $X_1<X_2<\cdots<X_l$ is the
set $\{\rank X_i\}_{1\le i\le l}$.  These numbers are the components
of the \emph{flag $f$-vector} of~$L$.  The \emph{flag $h$-vector} is
defined by
\begin{equation} \label{e:flagh}
h_J(L)=\sum_{K\subseteq J} (-1)^{\abs{J-K}} f_K(L),
\end{equation}
or, equivalently,
\begin{equation} \label{e:flagf}
f_J(L)=\sum_{K\subseteq J} h_K(L).
\end{equation}
Up to sign, $h_J(L)$ is the reduced Euler characteristic of the
rank-selected subposet $L_J$ of~$L$.  More precisely,
\begin{equation} \label{e:ec}
h_J(L)=(-1)^{\abs{J}-1}\tilde\chi(L_J).
\end{equation}
If the order complex of $L_J$ is homotopy equivalent to a wedge of
$(\abs{J}-1)$-spheres, then $h_J(L)$ is the number of spheres.

For the Boolean lattice, one can see from \eqref{e:flagf} that
$h_J(L)$ is equal to the descent number $\beta(J)$.  This is also a
special case of Proposition~\ref{p:beta} (Section~\ref{sub:flagh}).
The main result of this subsection, generalizing \eqref{e:descent}, is
the following proposition.

\begin{proposition}[Stanley, private communication] \label{p:stanley}
Let $L$ be a graded poset with $\zerohat,\onehat$, and let $n$ be its
rank.  Then
\begin{equation} \label{e:descent1}
d(L)=\sum_{\substack{J\subseteq[n-1]\\ J\in\J}} h_J(L).
\end{equation}
\end{proposition}

\begin{proof}
Let $d'(L)$ denote the right-hand side of \eqref{e:descent1} if
$L\ne0$, and let $d'(0)=1$.  It suffices to show that $d'$ satisfies
the recurrence \eqref{e:recurrence2}, i.e.,
\begin{equation} \label{e:recurrence4}
d'(L)=\mu_L(\zerohat,\onehat) + \sum_{X\in M} d'([\zerohat,X]),
\end{equation}
where $M$ is the set of elements of $L$ of rank $n-1$.  We may assume
$n\ge2$.  Group the terms on the right-hand side of \eqref{e:descent1}
in pairs, where $J\subseteq[n-2]$ is paired with $J_+=J\cup\{n-1\}$.
This leaves one term unpaired: If $n$ is even, we have
$[n-1]=[n-2]_+\in\J$ but $[n-2]\notin\J$, while the reverse is
true if $n$ is odd.  In both cases we obtain
\begin{equation}
d'(L)=(-1)^n h_{[n-1]}(L) +
\sum_{\substack{J\subseteq[n-2]\\ J\in\J}}\left(h_J(L)+h_{J_+}(L)\right). 
\end{equation}
Two simple observations now complete the proof of
\eqref{e:recurrence4}.  The first is that
\[
(-1)^n h_{[n-1]}(L) = \mu_L(\zerohat,\onehat)
\]
by \eqref{e:ec} with $J=[n-1]$.  The second observation is that
\[
h_J(L)+h_{J_+}(L)=\sum_{X\in M} h_J([\zerohat,X])
\]
for $J\subseteq[n-2]$.  This is proved by expanding both terms on the
left-hand side by \eqref{e:flagh}, noting that many terms cancel, and
applying the following fact to the remaining terms:
\[
f_{K_+}(L)=\sum_{X\in M} f_K([\zerohat,X])
\]
for $K\subseteq[n-2]$.
\end{proof}

Michelle Wachs [private communication] has pointed out that
Proposition~\ref{p:stanley} implies the following result about
$q$-derangement numbers, due to D\'esarm\'enien and Wachs
\cite[Section 7]{\desarw}:

\begin{corollary}
The $q$-derangement number $d_n(q)$ satisfies
\[
d_n(q)=\sum_{\pi\in E_n} q^{\inv(\pi)},
\]
where $E_n$ is the set of desarrangements in $S_n$ and $\inv(\pi)$ is
the number of inversions of $\pi$.
\end{corollary}

\begin{proof}
Take $L$ to be the subspace lattice of~$\fq^n$, so that $d(L)=d_n(q)$.
It is known \cite[Theorem 3.12.3]{\stanley} that
\[
h_J(L)=\sum_{\substack{\pi\in S_n\\ \des(\pi)=J}} q^{\inv(\pi)}.
\]
The corollary now follows at once from the proposition.
\end{proof}

\subsection{More on the flag $h$-vector} \label{sub:more}

Going back to the random walk on maximal chains for motivation, recall
that there is an eigenvalue $\lambda_X$ for each $X\in L$ (where $L$
is the lattice of flats of a matroid), with multiplicity
$m_X=d([X,\onehat])$.  We have just seen that this multiplicity is a
sum of certain components of the flag $h$-vector when $X=\zerohat$.
Is the same true of the other multiplicities?  This is a reasonable
question since
\[
\sum_{X\in L} m_X = f(L) = \sum_{J\subseteq[n-1]} h_J(L).
\]
One might naively hope to lump the terms on the right-hand-side in
such a way that each lump accounts for one~$m_X$.  This does not seem
to be the case; but what is true is that if we lump together all the
$m_X$ with $X$ of a given rank, then their sum is equal to the sum of
the $h_J$ for certain sets $J$.  This was observed by Swapneel Mahajan
[private communication].  It is of interest for the random walk in
case $L$ has the property that all flats of a given rank contain the
same number of rank~1 flats.  If, further, we take uniform weights on
the rank~1 flats, then we get one eigenvalue for each possible
rank~$r$, $0\le r\le n=\rank(L)$, the multiplicity being
\[
D_{n-r}(L) := \sum_{\rank(X)=r} d([X,\onehat]).
\]
(The subscript $n-r$ is a reminder that each interval $[X,\onehat]$ on the
right has rank $n-r$.)  Mahajan's result, then, is that $D_{n-r}(L)$
is a sum of certain values of the flag $h$-vector.  This is valid for
every graded poset with $\zerohat,\onehat$.  When $r=0$ it reduces to
Stanley's result from the previous section.

To state the result precisely, we associate to every set
$J\subseteq[n-1]$ a number $\gamma=\gamma(J)$, $0\le\gamma\le n$, as
follows.  Arrange the elements of $J$ in order, and consider the
initial run of consecutive integers; this has the form
$i,i+1,\dots,i+l-1$, where $l$ is the length of the run.  We allow the
case $J=\emptyset$, in which case we set $l=0$ and $i=n$.  Then $\gamma$
is defined by
\[
\gamma(J)=\begin{cases}
i&\text{if $l$ is even}\\
i-1&\text{if $l$ is odd.}
\end{cases}
\]
The result, then, is:

\begin{proposition}[Mahajan, private communication]
If $L$ is a graded poset with $\zerohat,\onehat$ and $n=\rank(L)$, then
\[
D_{n-r}(L)=\sum_{\gamma(J)=r} h_J(L).
\]
\end{proposition}

We omit the proof.  The starting point is to apply
Proposition~\ref{p:stanley}  to each of the posets $[X,\onehat]$.

For our random walk, the proposition says that the total multiplicity
of the eigenvalues contributed by the $X\in L$ of rank~$r$ is given by
the components of the flag $h$-vector with $\gamma(J)=r$.

Here are some special cases.
\begin{itemize}
\item
$r=0$:  We have $\gamma(J)=0$ if and only if the initial run in $J$ is
$1,\dots,l$ with $l$ odd, so the first omitted integer is even, as in
Proposition~\ref{p:stanley}.
\item
$r=n-1$:  There is no $J$ with $\gamma(J)=n-1$, so $D_1(L)=0$.  This is
consistent with the fact that $d=0$ for posets of rank~$1$.
\item
$r=n$:  The only $J$ with $\gamma(J)=n$ is $J=\emptyset$, so
$D_0(L)=h_\emptyset(L)=1$.  This is consistent with the fact that
$d=1$ for the trivial poset.
\end{itemize}

\bibliographystyle{hamsplain}
\bibliography{hyperplane}

\end{document}